\documentclass[11pt]{article}

\usepackage{amsmath}
\usepackage{amssymb}
\usepackage{amsthm}
\usepackage{mathtools}
\usepackage{bm}
\usepackage{graphicx}
\usepackage{hyperref}

\usepackage[margin=1in]{geometry}
\usepackage{amsmath} 
\usepackage{amssymb} 
\usepackage{amsbsy} 
\usepackage{amsthm} 
\usepackage{amsfonts}
\usepackage{thmtools}
\usepackage{thm-restate}
\usepackage{bbm}
\usepackage{url}
\usepackage{booktabs}
\usepackage{mathrsfs}
\usepackage{multirow}
\usepackage{chngpage}
\usepackage{subcaption}
\usepackage{float}
\usepackage{subfiles}
\usepackage[math]{cellspace}
\usepackage[font={small}]{caption}
\usepackage{tikz, pgfplots}
\usepackage{algpseudocode}
\usepackage{graphicx}
\usepackage{algorithm}
\usepackage{fullpage}
\usepackage{chngcntr}
\usepackage[square,sort,comma,numbers]{natbib}
\usepackage{soul}
\usepackage{mathtools}
\usepackage{array}
\usepackage{arydshln}
\usepackage[letterspace=-85]{microtype}
\usepackage{nicefrac}
\usepackage{hyperref}

\usetikzlibrary{shapes.geometric,arrows,positioning,calc}
\usetikzlibrary{patterns}
\usetikzlibrary{decorations.pathmorphing,backgrounds,fit,matrix,decorations.pathreplacing}

\newcommand*{\pfstart}{\begin{proof}}
\newcommand*{\pfend}{\end{proof}}

\tikzstyle{vertex} = [circle, minimum size=0.1cm, inner sep=0pt, draw=black, fill=black]
\tikzstyle{circ} = [circle, minimum width=0.5mm, inner sep=0pt,draw,fill]
\tikzstyle{hcirc} = [circle,minimum width=5mm, inner sep=0pt,draw]
\tikzstyle{bcirc} = [circle, minimum width=1.5mm, inner sep=0pt,draw,fill]
\tikzstyle{bhcirc} = [circle, minimum width=1.5mm, inner sep=0pt,draw, dotted ]
\tikzstyle{ept} = [circle,minimum width=0mm, inner sep=0pt, white]
\tikzstyle{txt} = [text width=1.3cm,draw,rounded corners=3pt]
\tikzstyle{ncirc} = [circle,draw=black, inner sep=1pt]

\DeclarePairedDelimiter\iprod{\langle}{\rangle}

\algrenewcommand\algorithmicindent{1em}

\theoremstyle{remark}

\usepackage{bm}
\theoremstyle{remark}
\newtheorem*{claim*}{Claim}
\theoremstyle{remark}
\newtheorem*{remark*}{Remark}
\theoremstyle{remark}
\newtheorem{remark}{Remark}

\theoremstyle{plain}

\theoremstyle{plain}

\theoremstyle{definition}

\theoremstyle{definition}
\newtheorem*{assumption*}{Assumption}

\theoremstyle{definition}
\newtheorem*{example*}{Example}

\renewcommand\thmcontinues[1]{Cont.}

\theoremstyle{definition}
\newtheorem{definition}{Definition}

\theoremstyle{plain}

\theoremstyle{plain}

\usepackage{tikz}
\usetikzlibrary{arrows.meta,calc,decorations.pathmorphing}

\providecommand{\pfend}{}
\renewcommand{\pfend}{\hfill\ensuremath{\square}\par}

\graphicspath{{Figures/}}  

\usepackage{adjustbox}
\usepackage{natbib} 

\usepackage{float}
\usepackage{placeins}

\usepackage{adjustbox}
\usepackage{caption}

\newcommand{\TABLE}[3]{%
  \caption{#1}%
  \vspace{2pt}
  \noindent
  \begin{adjustbox}{max width=\textwidth}
  #2%
  \end{adjustbox}%
  \par\vspace{2pt}
  \noindent{\footnotesize #3}%
}

\newcommand{\FIGURE}[3]{%
  \caption{#2}%
  \begin{center}
    #1%
  \end{center}
  \par\noindent{\footnotesize\emph{Note.} #3}%
}

\theoremstyle{definition}

\usepackage{amsthm}
\usepackage{enumitem}

\hypersetup{
    colorlinks=true,
    citecolor=black,
    linkcolor=black,
    filecolor=magenta,      
    urlcolor=cyan,
}

\pgfplotsset{compat=1.15}
\tikzset{snake it/.style={decorate, decoration={snake, amplitude=.3mm, segment length=2mm}}}

\allowdisplaybreaks

\begin{document}

\title{Scalable Finite Adaptability via Polyhedral Partition and Learning}
	\author{
     Zolykha Rezaei
    \thanks{
    Department of Industrial, Manufacturing \& Systems Engineering, Texas Tech University, Lubbock, TX 79409 (email: {\tt zrezaei@ttu.edu}).}
    \and
    Ningji Wei
    \thanks{
    Department of Industrial, Manufacturing \& Systems Engineering, Texas Tech University, Lubbock, TX 79409 (email: {\tt ningji.wei@ttu.edu}). This work was partially supported by AFOSR grant FA9550-26-1-B345 and NSF grant 2548574.}
    \and
    Eojin Han
    \thanks{
        Information Technology, Analytics and Operations Mendoza College of Business, University of Notre Dame, Notre Dame, IN 46556 (email: {\tt ehan3@nd.edu}).}
}
	\date{}
	\maketitle
\begin{abstract}
We study finite adaptability for decision-making under uncertainty, where a finite set of candidate decisions is prepared in advance and a response is selected after uncertainty is realized. Despite advances in exact formulations, scalability remains challenging due to combinatorial decision assignments and the joint optimization of uncertainty set partitions and associated decisions. We develop a framework that explicitly parameterizes polyhedral partitions, enabling direct optimization of partition boundaries and associated decisions. We establish convergence of such polyhedral policies to the fully adjustable benchmark under general risk measures, including a rate of $O(K^{-\alpha/n})$ under uniform $\alpha$-H\"older continuity, where $K$ is the number of regions and $n$ is the uncertainty dimension. 
Exploiting the explicit polyhedral structure, we derive two exact reformulations for robust finite adaptability, one based on master–subproblem cut generation and the other on piecewise dualization, and extend the framework to stochastic programming.
To improve scalability, we introduce an approximate-learn-parallel algorithm that integrates partition learning with parallel regional optimization while preserving robustness. Computational experiments on classical robust and stochastic testbeds demonstrate scalability to larger instances and competitive policy performance.
\par\vspace{0.7em}
\end{abstract}

\section{Introduction}\label{sec:Intro}

Decision-making under uncertainty arises naturally in many domains, including logistics, production planning, finance, and emergency response. A central challenge is balancing implementation readiness with response optimality. 
While ex ante approaches commit to a solution before uncertainty is resolved and ex post approaches defer decisions until uncertainty is revealed, contingency planning offers a middle ground by preparing a finite menu of candidate solutions in advance, enabling rapid adaptation to realized conditions while preserving implementation readiness.



This idea of contingency planning has been formalized in several related optimization frameworks. In robust optimization (RO), $K$-adaptability prepares $K$ candidate recourse decisions to approximate fully adaptive policies in two-stage or multistage problems \citep{bertsimas2010finite,hanasusanto2015k, subramanyam2020kadaptability}. Related min--max--min formulations focus on selecting among candidate responses without first-stage decisions \citep{buchheim2017minmaxmin,goerigk2020minmaxmin, kurtz2024approximation}. While these approaches typically fix $K$ in advance, adaptive partitioning methods dynamically refine the uncertainty set and associated decisions \citep{bertsimas2016adaptivepartitions,postek2016iterative}. In stochastic programming (SP), related min--$\mathbb{E}$--min formulations typically consider finite scenario sets \citep{buchheim2019k,malaguti2022k}, where partitioning amounts to scenario-to-policy assignment. 
Consequently, methods for robust and stochastic settings have largely evolved independently. Finite adaptability has also been studied in distributionally robust optimization (DRO), where contingent decisions are assigned to partitions of the uncertainty support \citep{hanasusanto2016k,han2023finiteadaptability}.

These formulations share a common structure: the uncertainty set is partitioned, implicitly or explicitly, into regions associated with individual decisions, which we collectively refer to as \emph{finite adaptability}. This paper focuses on the robust (min--max--min) and stochastic (min--$\mathbb{E}$--min) settings without first-stage decisions. Table~\ref{tab:finite_adaptability_lit} provides a detailed comparison with existing formulations.


Significant progress has been made in developing solution methods for finite adaptability, including mixed-integer reformulations \citep{hanasusanto2015k}, branch-and-bound algorithms \citep{subramanyam2020kadaptability}, and iterative partition refinement schemes \citep{bertsimas2016adaptivepartitions, postek2016iterative}. Nevertheless, scalability remains a persistent challenge due to three interrelated sources of complexity. First, assigning scenarios to candidate responses induces combinatorial complexity that grows rapidly with the number of scenarios and candidate decisions. Second, continuous uncertainty sets introduce semi-infinite constraints, requiring robustness over infinitely many realizations. Third, partition construction and response optimization are coupled, requiring the joint optimization of all $K$ responses and their associated regions. Consequently, exact methods can struggle beyond moderate problem sizes. For example, in capital budgeting problems with four-dimensional uncertainty and $30$ projects, existing exact approaches reach a two-hour time limit on most test instances with $K=4$ \citep{hanasusanto2015k,subramanyam2020kadaptability}.

These challenges motivate a geometric perspective on finite adaptability. While existing approaches often treat the selection among candidate decisions algebraically, we make the underlying partition of the uncertainty set explicit, enabling new formulations and scalable optimization strategies. Moreover, since finite adaptability itself approximates fully adjustable policies, an approximate solution with more candidate responses may outperform an exact solution restricted to fewer responses. Building on these observations, three contributions follow.


On the theoretical side, we establish convergence of polyhedral policies to the fully adjustable benchmark under general risk measures. In particular, we establish a convergence rate of $O(K^{-\alpha/n})$ under uniform $\alpha$-H\"older continuity, where $K$ denotes the number of partition regions and $n$ the uncertainty dimension. These guarantees extend to two-stage problems (Theorem~\ref{thm:convergence} and Corollary~\ref{cor:two_stage}).

On the formulation side, we develop a parametric polyhedral partition framework that enables direct optimization of partitions alongside their associated decisions. Exploiting the polyhedral structure, we develop a master--subproblem method and, to our knowledge, the first piecewise-dualization reformulation for robust finite adaptability, both yielding finite-dimensional formulations.

On the computational side, we develop a learning-based approach that combines partition learning with parallel optimization while preserving robustness. The learned partition can be further refined using the piecewise-dualization reformulation. Computational experiments on classical testbeds demonstrate scalability and competitive policy performance in both robust and stochastic settings.

\subsection{Literature Review}
\label{sec:lit}

\begin{table}[!tbp]
\TABLE
{\raggedright Overview of finite adaptability in RO, SP, and DRO. \label{tab:finite_adaptability_lit}}
{\small\setlength{\tabcolsep}{4pt}\renewcommand{\arraystretch}{1.2}\centering
\begin{adjustbox}{max width=\linewidth,center}
\begin{tabular}{lccccc}
\hline
Reference & Type & \(K\) & Model & Uncertainty & Recourse Variables \\
\hline
\citet{bertsimas2010finite} & RO & Fixed & Multistage  & General & General \\
\citet{hanasusanto2015k} & RO & Fixed & Two-stage & General & Binary \\
\citet{bertsimas2015design} & RO & Fixed & Multistage & General & General \\
\citet{bertsimas2016adaptivepartitions} & RO & Not fixed & Multistage & General & General \\
\citet{postek2016iterative} & RO & Not fixed & Multistage & General & General \\
\citet{buchheim2017minmaxmin} & RO & Fixed & Min--max--min & Objective & Binary \\
\citet{subramanyam2020kadaptability} & RO & Fixed & Two-stage & General & General \\
\citet{buchheim2019k} & SP & Fixed & Min--$\mathbb E$--min  & Objective & Binary \\
\citet{malaguti2022k} & SP & Fixed & Min--$\mathbb E$--min & General & Continuous  \\
\citet{donninik} & SP & Fixed & Two-stage & General & General  \\
\citet{hanasusanto2016k} & DRO & Fixed & Two-stage & General & Binary \\
\citet{han2023finiteadaptability} & DRO & Fixed & Two-stage & RHS & General \\
This paper & RO / SP & Not Fixed &Min--max--min / Min--$\mathbb E$--min & General & General \\
\hline
\end{tabular}
\end{adjustbox}}
{\begingroup\fontsize{8}{11}\selectfont\emph{Notes.}
The $K$ column indicates whether the number of candidate policies is fixed. ``General'' in the uncertainty column indicates uncertainty beyond the objective or right-hand side; in the recourse column, it indicates applicability to both integer and continuous variables.
\par\endgroup}
\end{table}

Adjustable RO allows decisions to adapt to realized uncertainty \citep{ben2004adjustable}, with decision-rules balancing policy performance and tractability. Common policy classes in RO and SP include affine \citep{bertsimas2010optimality,bertsimas2012power,kuhn2011primal,chen2008linear}, lifted affine \citep{georghiou2015generalized,georghiou2020primal}, and piecewise-affine policies \citep{ben2020tractable,han2025nonlinear,thoma2026note}, while dual decision rules approximate Lagrangian multipliers \citep{kuhn2011primal,daryalal2024lagrangian}. For a broader review, see \citet{yanikouglu2019survey}.

Among these decision-rule families, finite adaptability corresponds to piecewise-constant policies, making it particularly suitable for discrete decisions and practical implementation. Following its introduction by \citet{bertsimas2010finite}, \citet{bertsimas2011geometric} analyzed its geometric properties in multistage stochastic and adaptive optimization, while \citet{bertsimas2015design} demonstrated the effectiveness of restricted adaptability in multistage settings. For two-stage RO, \citet{hanasusanto2015k} developed mixed-integer formulations for binary recourse, and \citet{subramanyam2020kadaptability} proposed an exact partition-based branch-and-bound algorithm for mixed-integer problems.
More recently, \citet{julien2025machine} employed machine learning to guide node selection in $K$-adaptability
branch-and-bound.
In contrast to approaches with a fixed number of candidate decisions, \citet{bertsimas2016adaptivepartitions} and \citet{postek2016iterative} developed adaptive partition-refinement schemes. Related min--max--min RO considers candidate selection without first-stage decisions \citep{buchheim2017minmaxmin,buchheim2018complexity}, with subsequent extensions to broader settings \citep{goerigk2020minmaxmin,kurtz2024approximation}.


Finite adaptability has also been studied in SP and DRO. In SP, \citet{buchheim2019k} study min--$\mathbb E$--min problems under objective uncertainty, while \citet{malaguti2022k} analyze their complexity and develop exact algorithms. More generally, \citet{donninik} investigate two-stage stochastic problems with objective and constraint uncertainty, providing approximation guarantees and exact partition-based methods. In DRO, \citet{hanasusanto2016k} study binary recourse, while \citet{han2023finiteadaptability} develop structured partitions for two-stage problems. Table~\ref{tab:finite_adaptability_lit} summarizes the literature across optimization paradigms, model structures, and uncertainty settings.

The remainder of the paper is organized as follows. Section~\ref{sec:pre} establishes convergence properties of polyhedral policies. Section~\ref{sec:model} introduces the parametric partition framework and reformulations. Section~\ref{sec:learn} develops a scalable learning-based solution framework, which is extended to SP in Section~\ref{sec:extsp}. Section~\ref{sec:experiments} presents computational experiments, and Section~\ref{sec:conclusion} concludes. Throughout, bold lowercase and uppercase letters denote vectors and matrices. $\bm{1}$ is the all-ones vector, $\odot$ denotes element-wise multiplication, $\iprod{\cdot,\cdot}$ denotes the inner product, and $[q]=\{1,\ldots,q\}$.

\section{Convergence of Polyhedral Policies}
\label{sec:pre}
We consider a decision problem under uncertainty with solution space $\mathcal Y$, a convex and compact uncertainty set $\Xi$, and a cost function $f:\mathcal Y\times\Xi\to \mathbb R\cup\{\pm\infty\}$. We adopt extended-value conventions for infeasibility and unboundedness and assume that all problems with finite optimal values attain their optima. The \emph{fully adjustable problem} is
\begin{equation}
  \label{eq:adj}
  \min_{\pi(\cdot)\in\mathcal Y^\Xi}
  \mathcal R_\Xi\left(f\!\left(\pi(\bm\xi),\bm\xi\right)\right),
\end{equation}
where $\mathcal Y^\Xi$ denotes the space of policies
$\pi:\Xi\to\mathcal Y$, and the risk measure
$\mathcal R_\Xi$ \citep{shapiro2021lectures} aggregates
costs across uncertainty realizations. Examples include
the worst-case operator $\mathcal R_\Xi(u)
=\sup_{\bm\xi\in\Xi}u(\bm\xi)$ in RO,
the expectation $\mathcal R_\Xi(u)
=\mathbb E_{\bm\xi\sim\mathbb P}[u(\bm\xi)]$
in SP, and the worst-case expectation
over an ambiguity set in DRO.
The risk measure $\mathcal R_\Xi$ satisfies three standard properties:
\emph{monotonicity}, $u\leq v$ implies $\mathcal R_\Xi(u)\leq\mathcal R_\Xi(v)$;
\emph{normalization}, $\mathcal R_\Xi(0)=0$; and
\emph{translation equivariance},
$\mathcal R_\Xi(u+\alpha)=\mathcal R_\Xi(u)+\alpha$ for any constant $\alpha$.

By monotonicity of $\mathcal R_\Xi$, an optimal fully adjustable policy can be obtained through pointwise optimization, $\pi(\bm\xi)\in\arg\min_{\bm y\in\mathcal Y}f(\bm y,\bm\xi)$. However, fully adjustable policies are generally computationally intractable \citep{ben2004adjustable,ben2009robust}. Finite adaptability provides a tractable policy class by preparing $K$ candidate decisions in advance and selecting one after uncertainty is realized \citep{bertsimas2010finite}. It can be formulated as
\begin{equation}
  \label{eq:kadpt}
  \min_{\bm y \in \mathcal Y^{\mathcal K},\;
  \sigma(\cdot)\in\mathcal K^\Xi}
  \mathcal R_{\Xi}\left(
  f\!\left(\bm y_{\sigma(\bm\xi)},\bm\xi\right)
  \right),
\end{equation}
where $\mathcal K:=[K]$,
$\bm y=(\bm y_k)_{k\in\mathcal K}$ denotes the
$K$ candidate decisions, and $\sigma:\Xi\to\mathcal K$
is the selection policy.
Each selection policy $\sigma$ induces a partition
$\{\sigma^{-1}(k)\}_{k\in\mathcal K}$ of $\Xi$.
Conversely, a partition determines a selection policy
through a tie-breaking rule. Thus, \eqref{eq:kadpt}
is equivalently formulated as
\begin{equation}
  \label{eq:kpart}
  \min_{\substack{
  \bm y:=(\bm y_k)_{k \in \mathcal K} \in \mathcal Y^{\mathcal K} \\
  \mathcal P:=(\Xi_k)_{k \in \mathcal K} \in \Pi_K(\Xi)}}
  \mathcal R_{\Xi}\left(
  f\!\left(\bm y_{\sigma_{\mathcal P}(\bm\xi)},\bm\xi\right)
  \right),
\end{equation}
where $\Pi_K(\Xi)$ denotes the collection of
$K$-partitions of $\Xi$, and
$\sigma_{\mathcal P}(\bm\xi)
:=\min\{k\in\mathcal K:\bm\xi\in\Xi_k\}$
specifies the selection rule.
Since empty regions do not affect the optimization and boundary overlaps can be resolved by the tie-breaking rule $\sigma_{\mathcal P}(\cdot)$, such collections induce well-defined policies and can be refined into disjoint partitions of $\Xi$. We therefore adopt the following convention.
\begin{remark}\label{lem:partition}
Throughout the paper, we use the term \emph{partition}
to refer to a collection of sets whose union equals $\Xi$,
allowing empty parts and boundary overlaps.
\end{remark}

The partition formulation~\eqref{eq:kpart} motivates studying the structure of optimal partitions. Prior work has shown that optimal polyhedral partitions exist under objective uncertainty \citep{hanasusanto2015k}, while finite unions of polyhedra can characterize optimal partitions under right-hand-side uncertainty \citep{han2023finiteadaptability}. The following proposition extends these structural results to more general affine and piecewise-affine cost functions.

\begin{restatable}{proposition}{piecewise}
\label{prop:piecewise}
Suppose that $f(\bm y,\cdot)$ is affine in $\bm\xi$ for every $\bm y\in\mathcal Y$. Then \eqref{eq:kpart} admits an optimal partition $(\Xi_k)_{k\in\mathcal K}$ such that each $\Xi_k$ is the intersection of $\Xi$ with a polyhedron. More generally, if $f(\bm y,\cdot)$ is piecewise affine in $\bm\xi$ with finitely many polyhedral pieces for every $\bm y\in\mathcal Y$, then an optimal partition exists with each $\Xi_k$ given by the intersection of $\Xi$ with a finite union of polyhedra.
\end{restatable}

This proposition establishes the optimality of polyhedral partitions under affine objective uncertainty. With affine constraint uncertainty, an optimal partition may require refinement into $K'\geq K$ polyhedral regions. We call such policies \emph{polyhedral policies}.
For more general settings, the following theorem establishes readily verifiable conditions under which polyhedral policies converge to the fully adjustable benchmark as $K\to\infty$. For any policy $\pi:\Xi\to\mathcal Y$, let $\phi_\pi(\bm\xi):=f(\pi(\bm\xi),\bm\xi)$ denote its value function, with $\mathcal R_\Xi(\phi_\pi)$ representing its aggregated value.

\begin{restatable}{theorem}{polyconv}
\label{thm:convergence}
Let $\Xi\subseteq\mathbb R^n$ be compact, and consider
an optimal fully adjustable policy
\[
\pi(\bm\xi)\in\arg\min_{\bm y\in\mathcal Y_0}
\left\{g(\bm y,\bm\xi):
h_j(\bm y,\bm\xi)\leq0,\ j\in[m]\right\}.
\]
Suppose $\phi_\pi$ is finite on $\Xi$,
$\mathcal Y_0$ is convex,
$g(\cdot,\bm\xi)$ and $h_j(\cdot,\bm\xi)$
are convex for all $\bm\xi\in\Xi$ and $j\in[m]$,
and every $\bm\xi\in\Xi$ admits a strictly feasible
$\bar{\bm y}_{\bm\xi}\in\mathcal Y_0$.
Define the family of functions
\(
\mathcal F:=
\{g(\bm y,\cdot)\mid\bm y\in\mathcal Y_0\}
\cup
\{h_j(\bm y,\cdot)\mid
\bm y\in\mathcal Y_0,\ j\in[m]\}.
\)
Under each of the following regularity conditions,
there exists a sequence of polyhedral policies
$(\pi_K)_{K\in\mathbb N}$ satisfying
the corresponding convergence guarantee:
\begin{enumerate}[label=(\roman*)]
    \item If every function in $\mathcal F$ is
    continuous on $\Xi$, then
    $\phi_{\pi_K}\to\phi_\pi$ pointwise.

    \item If, in addition to (i), $\mathcal Y_0$
    is compact and $\mathcal F$ is equicontinuous
    on $\Xi$, then $\phi_{\pi_K}\to\phi_\pi$
    uniformly and
    $\mathcal R_\Xi(\phi_{\pi_K})
    \to\mathcal R_\Xi(\phi_\pi)$.

    \item If, in addition to (ii), $\mathcal F$
    satisfies a uniform H\"older condition
    on $\Xi$ with exponent $\alpha\in(0,1]$,
    then both convergence results hold at rate
    $O(K^{-\alpha/n})$.
\end{enumerate}
\end{restatable}

When $\mathcal Y_0$ is compact, joint continuity of $g$ and $h_j$ on $\mathcal Y_0\times\Xi$ implies condition (ii). Existing convergence results for finite adaptability primarily address linear two-stage RO and rely on continuity conditions that can be difficult to verify~\citep{bertsimas2010finite,kedadsidhoum2025}. Theorem~\ref{thm:convergence} extends these results to convex settings under a general risk measure $\mathcal R_\Xi$, with verifiable conditions based on convexity, continuity, and strict feasibility. It also establishes an explicit approximation bound, yielding a convergence rate of $O(K^{-\alpha/n})$ under uniform H\"older continuity. The following corollary extends these guarantees to two-stage problems.

\begin{restatable}{corollary}{polytwostage}
\label{cor:two_stage}
Consider the two-stage problem \( z^\star=\inf_{\bm x\in\mathcal X} \{c(\bm x)+\mathcal R_\Xi(\phi_{\bm x})\}, \) where $\phi_{\bm x}$ denotes the optimal second-stage value function given $\bm x$. Suppose an optimal first-stage solution $\bm x^\star$ exists and its corresponding second-stage problem satisfies the common assumptions of Theorem~\ref{thm:convergence}. Then, under each of the regularity conditions (i)--(iii) of Theorem~\ref{thm:convergence}, there exists a sequence of polyhedral recourse policies satisfying the corresponding convergence guarantees. In particular, the optimal two-stage values $z_K$ satisfy $z_K\to z^\star$ under (ii), and $z_K-z^\star=O(K^{-\alpha/n})$ under (iii).
\end{restatable}



If the first-stage optimum is unattained, convergence of $z_K$ still holds along a minimizing sequence satisfying the corresponding recourse assumptions.
Motivated by these convergence guarantees, the remainder of the paper develops solution methods for polyhedral policies. 
We focus on settings where $\Xi$ is polyhedral and the objective and constraints are affine in $\bm\xi$. For discrete $\Xi$, partitioning a polyhedral outer approximation (e.g., its convex hull or LP relaxation) induces a partition of $\Xi$, extending the framework to discrete uncertainty.



\section{Finite Adaptability with Parametric Polyhedral Partition}\label{sec:model}

In this section, we study the following \emph{robust} \(K\)-adaptability problem in linear settings. Its extension to SP will be discussed in Section~\ref{sec:extsp}.
\begin{equation}
  \label{eq:linear_kpart}
  \begin{aligned}
  \min_{\substack{\bm y := (\bm y_k)_{k \in \mathcal K} \in \mathcal Y^{\mathcal K} \\
  \mathcal P := (\Xi_k)_{k \in \mathcal K} \in \Pi_K(\Xi)\\z \in \mathbb{R} }}&~z\\
  \text{s.t.} \quad &~
  \iprod{\bm{Q}(\bm{y}_k), \bm{\xi}} + \bm{q}(\bm{y}_k) \leq z, \quad \forall k \in \mathcal{K},~\bm \xi \in \Xi_k\\
              &~\iprod{\bm Q_j (\bm y_{k}), \bm\xi}
  + \bm q_j (\bm y_{k}) \leq 0,
  \quad \forall k \in \mathcal{K},~j \in J,~ \bm \xi \in \Xi_k. \\
  \end{aligned}
\end{equation}
Here, $J$ indexes the constraints, and the uncertainty set $\Xi:=\{\bm\xi\in\mathbb R^n\mid\bm B\bm\xi\leq\bm b\}$ is a nonempty compact polyhedron. The mappings $\bm Q,\bm Q_j:\mathcal Y\to\mathbb R^n$ and $q,q_j:\mathcal Y\to\mathbb R$ define the objective and constraint coefficients. Problem~\eqref{eq:linear_kpart} jointly optimizes a $K$-partition $\mathcal P=(\Xi_k)_{k\in\mathcal K}$ and decisions $\bm y=(\bm y_k)_{k\in\mathcal K}$, each robust over its corresponding region. We address this through the following \emph{parametric polyhedral partition} framework.

\subsection{Parametric Polyhedral Partition}
To ensure $\mathcal P=(\Xi_k)_{k \in \mathcal K}$ forms a partition of $\Xi$ with a finite number of parameters to optimize, we begin with the following definition.
\begin{definition}[Parametric Polyhedral Partition]
\label{def:ppp}
Let \(\mathbb{R}^\ell\) be a target space equipped with a predefined partition \(\mathcal L=(\mathcal L_k)_{k\in\mathcal K}\), where each region is defined as \(\mathcal L_k := \{\bm x \in \mathbb R^\ell \mid \bm B_k \bm x \le \bm b_k\}\), and let 
$\bm \xi \mapsto \bm T_{\bm \theta} \bm\xi + \bm{\tau}_{\bm \theta}$
be a parametric affine mapping with $\bm \theta \in \Theta \subseteq \mathbb R^d$.
We refer to the tuple \((\mathcal L, \bm T, \bm \tau)\) as a \emph{partition scheme}. Then, the sets \(\mathcal P_{\bm \theta}=(\Xi_k^{\bm \theta})_{k\in\mathcal K}\), with 
$\Xi_k^{\bm \theta}:=\{\bm \xi \in \Xi \mid \bm B_k (\bm T_{\bm \theta}\bm \xi + \bm \tau_{\bm\theta}) \le \bm b_k\},$
is called the \emph{parametric polyhedral partition} induced by \(\bm \theta\) under the scheme \((\mathcal L, \bm T, \bm \tau)\).
\end{definition}

  Although $(\bm{B}_k, \bm{b}_k)$ may appear redundant, as their role can be absorbed into $(\bm{T}_{\bm{\theta}}, \bm{\tau}_{\bm{\theta}})$, enforcing partition constraints directly on $(\bm{T}_{\bm{\theta}}, \bm{\tau}_{\bm{\theta}})$ is difficult in general. Fixing a partition through an explicit specification of $(\bm{B}_k, \bm{b}_k)$ addresses this. Under a given scheme $(\mathcal L, \bm T, \bm \tau)$, \eqref{eq:linear_kpart} can be rewritten as
\begin{subequations}
  \label{eq:ppp_linear_kpart}
  \begin{align}
  \min_{\substack{\bm y := (\bm y_k)_{k \in \mathcal K} \in \mathcal Y^{\mathcal K} \\
  \bm \theta \in \Theta, z \in \mathbb{R}}} &~z\\
  \text{s.t.} \quad &~
  \iprod{\bm{Q}(\bm{y}_k), \bm{\xi}} + \bm{q}(\bm{y}_k) \leq z, \quad \forall k \in \mathcal{K},~\bm \xi \in \Xi_k^{\bm\theta}\label{eq:ppp_linear_kpart01}\\
              &~\iprod{\bm Q_j (\bm y_{k}), \bm\xi}
  + \bm q_j (\bm y_{k}) \leq 0,
  \quad \forall k \in \mathcal{K},~j \in J,~ \bm \xi \in \Xi_k^{\bm \theta}, \label{eq:ppp_linear_kpart02}
  \end{align}
\end{subequations}
where each region $\Xi_k^{\bm \theta}
:= \{\bm \xi \in \mathbb R^n \mid \bm B \bm \xi \le \bm b,\ 
\bm B_k \bm T_{\bm \theta} \bm \xi \le \bm b_k - \bm B_k \bm \tau_{\bm \theta}\}$
is a parametric polytope. Different choices of partition schemes give rise to distinct formulations and, consequently, distinct classes of policies. In the following subsection, we introduce three specific partition schemes.

\subsection{Three Partition Schemes}
This subsection introduces three partition schemes: one orthant-based and two tree-based.
These approaches are also closely related to classical learning methods such as support vector machines (SVMs) and decision trees (DTs), enabling learning-based algorithms developed later in Section~\ref{sec:learn}.


\subsubsection{Orthant-Based Partition Scheme}

Figure~\ref{fig:orthant} illustrates the orthant-based partition scheme, where an affine mapping $\bm P\bm\xi+\bm p$ embeds $\Xi\subseteq\mathbb R^2$ into $\mathbb R^3$. The orthants of the target space induce a partition of $\Xi$ parameterized by $\bm\theta=(\bm P,\bm p)$, as formalized below.

\begin{figure}[!t]
\FIGURE
{%
\centering
\begin{subfigure}[t]{\textwidth}
  \centering
\resizebox{0.8\linewidth}{!}{
\begingroup%
\input{Figures/figure_style}%
\def\orthXHalfExtent{2.5}
\def\orthYHalfExtent{3.2}
\def\orthZBelow{2.2}
\def\orthZAbove{2.6}
\def\orthMiddleScale{.54}%
\def\orthViewSlope{.30}%
\newcommand{\orthSetupPlanes}{%
\pgfmathsetlengthmacro{\orthViewRise}{\orthViewSlope*1cm}
\pgfmathsetlengthmacro{\orthSurfaceRise}{(.55+\orthViewSlope)*1cm}
\pgfmathsetlengthmacro{\orthSurfaceOffset}{(-.4125-2.2*\orthViewSlope)*1cm}
\pgfmathsetmacro{\orthGUmin}{-\orthXHalfExtent}
\pgfmathsetmacro{\orthGUmax}{\orthXHalfExtent}
\pgfmathsetmacro{\orthGVmin}{-\orthYHalfExtent}
\pgfmathsetmacro{\orthGVmax}{\orthYHalfExtent}
\pgfmathsetmacro{\orthUmin}{-\orthXHalfExtent}
\pgfmathsetmacro{\orthUmax}{\orthXHalfExtent}
\pgfmathsetmacro{\orthVmin}{-\orthYHalfExtent}
\pgfmathsetmacro{\orthVmax}{\orthYHalfExtent}
\pgfmathsetmacro{\orthWmin}{-\orthZBelow}
\pgfmathsetmacro{\orthWmax}{\orthZAbove}
\pgfmathsetmacro{\orthVCutLo}{max(\orthVmin,\orthWmin/.55-1.45)}
\pgfmathsetmacro{\orthVCutHi}{min(\orthVmax,\orthWmax/.55-1.45)}
\pgfmathsetmacro{\orthUCutLo}{max(\orthUmin,\orthWmin/.55-1.45)}
\pgfmathsetmacro{\orthUCutHi}{min(\orthUmax,\orthWmax/.55-1.45)}
\pgfmathsetmacro{\orthGCutLo}{max(\orthGUmin,-1.45-\orthGVmax)}
\pgfmathsetmacro{\orthGCutHi}{min(\orthGUmax,-1.45-\orthGVmin)}
\pgfmathsetmacro{\orthVSeamLo}{max(\orthGVmin,\orthVmin)}
\pgfmathsetmacro{\orthVSeamHi}{min(\orthGVmax,\orthVmax)}
\pgfmathsetmacro{\orthUSeamLo}{max(\orthGUmin,\orthUmin)}
\pgfmathsetmacro{\orthUSeamHi}{min(\orthGUmax,\orthUmax)}
\pgfmathsetmacro{\orthScreenLo}{.6*min(\orthGUmin-.65*\orthGVmax,-.65*\orthVmax,\orthUmin)}
\pgfmathsetmacro{\orthScreenHi}{.6*max(\orthGUmax-.65*\orthGVmin,-.65*\orthVmin,\orthUmax)}
\pgfmathsetmacro{\orthMiddleX}{4.1-\orthScreenLo}
\pgfmathsetmacro{\orthMiddleY}{1.4-.5*\orthMiddleScale*(\orthWmax+\orthWmin)}
\pgfmathsetmacro{\orthArrowX}{\orthMiddleX+\orthScreenHi+.65}
\pgfmathsetmacro{\orthRightX}{\orthArrowX+1.4}
\pgfmathsetlengthmacro{\orthClipXmin}{(\orthScreenLo/.6-1)*1cm}
\pgfmathsetlengthmacro{\orthClipXmax}{(\orthScreenHi/.6+1)*1cm}
\pgfmathsetlengthmacro{\orthClipYmin}{(\orthWmin+\orthViewSlope*min(\orthUmin,\orthVmin)-1)*1cm}
\pgfmathsetlengthmacro{\orthClipYmax}{(\orthWmax+\orthViewSlope*max(\orthUmax,\orthVmax)+1)*1cm}
}%
\definecolor{orthSurfaceGray}{gray}{0.88}%
\newcommand{\orthVertices}{%
  \coordinate (A) at (0,0);
  \coordinate (B) at (1.75,0);
  \coordinate (C) at (.83,2.79);
  \coordinate (D) at (0,2.79);
  \coordinate (E) at (.8,0);
  \coordinate (F) at (.8,2.79);
  \coordinate (G) at (0,1.4);
  \coordinate (H) at (1.2883512545,1.4);
  \coordinate (I) at (.8,1.4);
  \coordinate (J) at (.75,0);
  \coordinate (K) at (0,.75);
}%
\newcommand{\orthFills}{%
  \fill[orthSurfaceGray] (A)--(B)--(C)--(D)--cycle;
}%
\newcommand{\orthRegions}{%
  \orthFills
  \draw[cut] (E)--(F) (G)--(H) (J)--(K);
  \draw[boundary] (A)--(B)--(C)--(D)--cycle;
}%
\newcommand{\orthLabels}{%
  \node[region] at (.24,.18) {$\Xi_1$};
  \node[region] at (.36,1.02) {$\Xi_2$};
  \node[region] at (1.13,.58) {$\Xi_3$};
  \node[region] at (.35,2.12) {$\Xi_4$};
  \node[region] at (1.015,1.68) {$\Xi_5$};
}%
\begin{tikzpicture}[
  >=Stealth,
  axis/.style={->,line width=.5pt},
  boundary/.style={line width=.7pt,black},
  hidden boundary/.style={line width=.55pt,black!70,dash pattern=on 2.3pt off 1.6pt},
  cut/.style={densely dotted,line width=.45pt,black},
  visible cut/.style={line width=.45pt,black},
  cut extension/.style={densely dotted,line width=.35pt,black!55},
  plane/.style={fill=black!6,fill opacity=.4,draw=black!24,line width=.3pt},
  plane seam/.style={black!30,line width=.3pt},
  region/.style={font=\faFigureDetailFont,inner sep=.5pt},
  surface/.style={shift={(.11cm,\orthSurfaceOffset)},x={(1cm,\orthSurfaceRise)},y={(-.65cm,\orthSurfaceRise)}},
  every node/.style={font=\faFigureFont,text=black,outer sep=2pt}
]
\orthSetupPlanes
\begin{scope}
  \orthVertices
  \orthFills
  \draw[boundary] (A)--(B)--(C)--(D)--cycle;
  \draw[axis] (0,0)--(2.05,0) node[right] {$\xi_1$};
  \draw[axis] (0,0)--(0,3.12) node[above] {$\xi_2$};
  \node at (.55,1.35) {$\Xi$};
\end{scope}
\draw[->,line width=.7pt] (2.65,1.45)--(3.45,1.45)
  node[midway,above] {$\bm P\bm\xi+\bm p$};
\begin{scope}[shift={(\orthMiddleX,\orthMiddleY)},scale=\orthMiddleScale,
  x={(1cm,\orthViewRise)},y={(-.65cm,\orthViewRise)},z={(0cm,1cm)}]
  \path[plane,fill=black!10] (\orthGUmin,\orthGVmin,0)--(\orthGUmax,\orthGVmin,0)
    --(\orthGUmax,\orthGVmax,0)--(\orthGUmin,\orthGVmax,0)--cycle;
  \path[plane] (0,\orthVmin,\orthWmin)--(0,\orthVmax,\orthWmin)
    --(0,\orthVmax,\orthWmax)--(0,\orthVmin,\orthWmax)--cycle;
  \path[plane] (\orthUmin,0,\orthWmin)--(\orthUmax,0,\orthWmin)
    --(\orthUmax,0,\orthWmax)--(\orthUmin,0,\orthWmax)--cycle;
  \draw[plane seam]
    (0,0,\orthWmin)--(0,0,\orthWmax)
    (0,\orthVSeamLo,0)--(0,\orthVSeamHi,0)
    (\orthUSeamLo,0,0)--(\orthUSeamHi,0,0);
  \begin{scope}[surface]
    \orthVertices
    \orthFills
  \end{scope}
  \begin{scope}[even odd rule,overlay]
  \clip (\orthClipXmin,\orthClipYmin) rectangle (\orthClipXmax,\orthClipYmax)
    (-.8,-1.4,-.4125)--(.95,-1.4,.55)--(.03,1.39,1.5785)
    --(-.8,1.39,1.122)--cycle;
  \draw[black!24,line width=.3pt]
    (0,\orthVmin,\orthWmin)--(0,\orthVmax,\orthWmin)
      --(0,\orthVmax,\orthWmax)--(0,\orthVmin,\orthWmax)--cycle
    (\orthUmin,0,\orthWmin)--(\orthUmax,0,\orthWmin)
      --(\orthUmax,0,\orthWmax)--(\orthUmin,0,\orthWmax)--cycle
    (\orthGUmin,\orthGVmin,0)--(\orthGUmax,\orthGVmin,0)
      --(\orthGUmax,\orthGVmax,0)--(\orthGUmin,\orthGVmax,0)--cycle;
  \draw[plane seam]
    (0,0,.7975)--(0,0,\orthWmax)
    (0,\orthVSeamLo,0)--(0,-1.45,0)
    (\orthUSeamLo,0,0)--(-1.45,0,0);
  \end{scope}
  \draw[cut extension]
    (0,\orthVCutLo,{.55*(\orthVCutLo+1.45)})--(0,-1.4,.0275)
    (0,1.39,1.562)--(0,\orthVCutHi,{.55*(\orthVCutHi+1.45)})
    (\orthUCutLo,0,{.55*(\orthUCutLo+1.45)})--(-.8,0,.3575)
    (.4883512545,0,1.06609319)--(\orthUCutHi,0,{.55*(\orthUCutHi+1.45)})
    (\orthGCutHi,{-1.45-\orthGCutHi},0)--(-.05,-1.4,0)
    (-.8,-.65,0)--(\orthGCutLo,{-1.45-\orthGCutLo},0);
  \begin{scope}[surface]
    \orthVertices
    \draw[visible cut] (E)--(I)--(G) (J)--(K);
    \draw[cut] (I)--(F) (I)--(H);
    \draw[hidden boundary] (E)--(B)--(C)--(D)--(G) (K)--(A)--(J);
    \draw[boundary] (J)--(E) (G)--(K);
    \node[region] at (.23,.21) {$1$};
    \node[region] at (.36,1.02) {$2$};
    \node[region] at (1.13,.58) {$3$};
    \node[region] at (.35,2.12) {$4$};
    \node[region] at (1.00,1.72) {$5$};
  \end{scope}
\end{scope}
\draw[->,line width=.7pt] (\orthArrowX,1.45)--({\orthArrowX+.8},1.45)
  node[midway,above] {pull back};
\begin{scope}[shift={(\orthRightX,0)}]
  \orthVertices
  \orthRegions
  \orthLabels
  \draw[axis] (0,0)--(2.05,0) node[right] {$\xi_1$};
  \draw[axis] (0,0)--(0,3.12) node[above] {$\xi_2$};
\end{scope}
\faNormalizeFigureWidth
\end{tikzpicture}%
\endgroup%
}%
\end{subfigure}
}
{\raggedright Illustration of the orthant-based partition scheme.\label{fig:orthant}}
{An uncertainty set $\Xi\subseteq \mathbb{R}^2$ is embedded into the target space $\mathbb{R}^3$ via $\bm P\bm\xi+\bm p$. The orthants of $\mathbb{R}^3$ partition the image $\bm P\Xi+\bm p$, and the inverse mapping induces a parametric partition of $\Xi$ into five regions, separated by the preimages of the orthant boundaries.}
\end{figure}

\begin{definition}[Orthant-Based Partition Scheme]
  \label{defi:orthant}
  Let the target space be $\mathbb R^\ell$, and let $\mathcal K := [2^\ell]$ index the $2^\ell$ orthants, each identified by a sign vector $\bm s_k \in \{1,-1\}^\ell$.  
The associated orthant-based partition $\mathcal L = (\mathcal L_k)_{k \in \mathcal K}$ is defined by
$\mathcal L_{k} := \{ \bm x \in \mathbb R^\ell \mid \bm s_k \odot \bm x \ge \bm 0 \}$. 
Let $\bm \theta = (\bm P, \bm p)$, the associated affine mapping is defined by $\bm T_{\bm\theta}\bm\xi + \bm \tau_{\bm\theta} := \bm P \bm\xi + \bm p$.
The induced parametric polyhedral partition is $\Xi^{\bm\theta}_{k} := \{ \bm\xi \in \Xi \mid \bm s_k \odot (\bm P \bm\xi + \bm p) \ge \bm 0 \}$ for each $k \in \mathcal K$.
\end{definition}

As illustrated in Figure~\ref{fig:orthant}, each row of $\bm P$, together with the corresponding entry of $\bm p$, defines a separating hyperplane in $\Xi$, yielding flexible partitions despite the fixed orthants in the target space. This mechanism is closely related to SVM-based hyperplane partitioning.

\subsubsection{Decision Tree Partition Scheme}
Inspired by decision trees in machine learning, this scheme constructs hierarchical polyhedral partitions of $\Xi$ through coordinate-wise splits with optimizable thresholds, as defined below.

\begin{figure}[!t]
\FIGURE
{%
\centering
\begin{subfigure}[b]{\textwidth}
  \centering
  \resizebox{0.9\linewidth}{!}{
\begingroup%
\input{Figures/figure_style}%
\def\dtVerticalPadding{5pt}%
\def\dtLevelSpacing{1.15}
\definecolor{dtSurfaceGray}{gray}{0.88}%
\newcommand{\dtCubeFaces}{%
  \fill[black!6] (0,0,1)--(1,0,1)--(1,1,1)--(0,1,1)--cycle;
  \fill[dtSurfaceGray] (0,0,0)--(1,0,0)--(1,0,1)--(0,0,1)--cycle;
  \fill[black!18] (1,0,0)--(1,1,0)--(1,1,1)--(1,0,1)--cycle;
}%
\newcommand{\dtCubeEdges}{%
  \draw[cube hidden] (0,0,0)--(0,1,0)--(1,1,0) (0,1,0)--(0,1,1);
  \draw[cube edge] (0,0,0)--(1,0,0)--(1,1,0)--(1,1,1)
    --(0,1,1)--(0,0,1)--cycle
    (0,0,1)--(1,0,1)--(1,1,1) (1,0,1)--(1,0,0);
}%
\newcommand{\dtLeaf}[5]{%
  \node[leaf] (#1) at (#2) {$\begin{gathered}
    \bm s^{U}_{#3}=[#4]\\[-1pt]
    \bm s^{L}_{#3}=[#5]
  \end{gathered}$};
}%
\begin{tikzpicture}[
  >=Stealth,
  every node/.style={font=\faFigureFont,text=black,outer sep=1.5pt},
  split/.style={draw=black!85,line width=.55pt,fill=white,
    rounded corners=1pt,minimum width=.76cm,minimum height=.44cm,inner sep=2pt},
  leaf/.style={font=\faFigureDetailFont,inner sep=2pt,align=center,fill=white},
  branch/.style={->,line width=.5pt,black!85},
  test/.style={font=\faFigureDetailFont,fill=white,inner sep=1pt,outer sep=1pt,yshift=2pt},
  threshold/.style={font=\faFigureDetailFont},
  cube edge/.style={line width=.65pt,black!90},
  cube hidden/.style={densely dotted,line width=.45pt,black!70},
  cube cut/.style={line width=.45pt,black!80},
  cut hidden/.style={densely dotted,line width=.45pt,black!75},
  cube view/.style={x={(1.10cm,-.28cm)},y={(.70cm,.45cm)},z={(0cm,1.25cm)}}
]
\begin{scope}[shift={(0,1.1)},cube view]
  \dtCubeFaces
  \dtCubeEdges
  \node at (.50,0,.64) {$\Xi$};
\end{scope}
\draw[->,line width=.7pt] (2.25,2.05)--(2.95,2.05);
\begin{scope}[shift={(0,{1.7*(1-\dtLevelSpacing)})},yscale=\dtLevelSpacing]
\node[split] (root) at (7.5,3.5) {$\xi_1$};
\node[split] (left) at (5.65,2.65) {$\xi_2$};
\node[split] (right) at (9.35,2.65) {$\xi_3$};
\node[split] (leftlow) at (4.5,1.6) {$\xi_3$};
\node[split] (righthigh) at (10.5,1.6) {$\xi_1$};
\dtLeaf{leaf1}{3.2,-.10}{1}{-1,-1,-1}{0,0,0}
\dtLeaf{leaf2}{5.8,-.10}{2}{-1,-1,0}{0,0,1}
\dtLeaf{leaf3}{6.3,1.55}{3}{-1,0,0}{0,1,0}
\dtLeaf{leaf4}{8.7,1.55}{4}{0,0,-1}{1,0,0}
\dtLeaf{leaf5}{9.3,-.10}{5}{-1,0,0}{1,0,1}
\dtLeaf{leaf6}{11.7,-.10}{6}{0,0,0}{1,0,1}
\draw[branch] (root.south)--(left.north)
  node[test,pos=.60,above] {$\xi_1\leq p_1$};
\draw[branch] (root.south)--(right.north)
  node[test,pos=.60,above] {$\xi_1\geq p_1$};
\draw[branch] (left.south)--(leftlow.north)
  node[test,pos=.48,left] {$\xi_2\leq p_2$};
\draw[branch] (left.south)--(leaf3.north)
  node[test,pos=.65,right,yshift=4pt] {$\xi_2\geq p_2$};
\draw[branch] (right.south)--(leaf4.north)
  node[test,pos=.65,left,yshift=4pt] {$\xi_3\leq p_3$};
\draw[branch] (right.south)--(righthigh.north)
  node[test,pos=.48,right] {$\xi_3\geq p_3$};
\draw[branch] (leftlow.south)--(leaf1.north)
  node[test,pos=.55,left] {$\xi_3\leq p'_3$};
\draw[branch] (leftlow.south)--(leaf2.north)
  node[test,pos=.55,right] {$\xi_3\geq p'_3$};
\draw[branch] (righthigh.south)--(leaf5.north)
  node[test,pos=.55,left] {$\xi_1\leq p'_1$};
\draw[branch] (righthigh.south)--(leaf6.north)
  node[test,pos=.55,right] {$\xi_1\geq p'_1$};
\end{scope}
\draw[->,line width=.7pt] (12.35,2.05)--(13.05,2.05);
\begin{scope}[shift={(13.5,1.1)},cube view]
  \dtCubeFaces
  \draw[cut hidden] (.38,0,0)--(.38,1,0)--(.38,1,1);
  \draw[cut hidden] (0,.55,1)--(0,.55,0)--(.38,.55,0)--(.38,.55,1);
  \draw[cut hidden] (0,0,.55)--(0,.55,.55)--(.38,.55,.55)--(.38,0,.55);
  \draw[cut hidden] (.38,0,.38)--(.38,1,.38)--(1,1,.38);
  \draw[cut hidden] (.72,0,.38)--(.72,1,.38)--(.72,1,1);
  \draw[cube cut]
    (.38,0,0)--(.38,0,1)--(.38,1,1)
    (0,.55,1)--(.38,.55,1)
    (0,0,.55)--(.38,0,.55)
    (.38,0,.38)--(1,0,.38)--(1,1,.38)
    (.72,0,.38)--(.72,0,1)--(.72,1,1);
  \dtCubeEdges
  \node[threshold,below] at (.38,0,-.03) {$p_1$};
  \node[threshold,above] at (.72,1,1.06) {$p'_1$};
  \node[threshold,above] at (0,.55,1.03) {$p_2$};
  \node[threshold,left] at (-.03,0,.55) {$p'_3$};
  \node[threshold,right] at (1.03,1,.38) {$p_3$};
\end{scope}
\path ([yshift=\dtVerticalPadding]current bounding box.north west)
  rectangle (current bounding box.south east);
\faNormalizeFigureWidth
\end{tikzpicture}%
\endgroup%
}%
\end{subfigure}
}
{\raggedright Illustration of the decision tree partition scheme.\label{fig:dt}}
{Given the domain $\Xi \subseteq \mathbb{R}^3$, a fixed decision tree routes each realization $\bm\xi$ via coordinate-wise threshold tests, producing a leaf $k \in \mathcal K$ with sign label $\bm s_k:=(\bm s^U_k, \bm s^L_k)$. This label encodes the root-to-leaf decisions and yields a polyhedral description of each leaf as $\Xi^{\bm\theta}_k=\{\bm \xi \in \Xi \mid \bm s^U_k \odot (\bm \xi - \bm p^U) \geq 0, \bm s^L_k \odot (\bm \xi - \bm p^L) \geq 0 \}$. Varying $\bm \theta$ induces different partitions on $\Xi$.
}
\end{figure}

\begin{definition}[Decision Tree Partition Scheme]
\label{defi:tree}
Given a tree with leaf set $\mathcal K$, let
$\bm s_k=(\bm s_k^U,\bm s_k^L)$ encode the
root-to-leaf splits for each $k\in\mathcal K$.
Specifically, $(\bm s_k^U)_i=-1$ if $\xi_i$
is bounded above along the path to $k$, and $0$
otherwise; similarly, $(\bm s_k^L)_i=1$ if
$\xi_i$ is bounded below, and $0$ otherwise.
Let $\ell=2n$ and define the target-space partition
$\mathcal L=(\mathcal L_k)_{k\in\mathcal K}$ by
$\mathcal L_k:=\{\bm x\in\mathbb R^\ell
\mid\bm s_k\odot\bm x\geq\bm0\}$.
With $\bm\theta=(\bm p^U,\bm p^L)$, define
the affine mapping
$\bm T_{\bm\theta}\bm\xi+\bm\tau_{\bm\theta}
:=(\bm\xi-\bm p^U,\bm\xi-\bm p^L)$.
The induced region is
$\Xi_k^{\bm\theta}
:=\{\bm\xi\in\Xi\mid
\bm s_k\odot(\bm\xi-\bm p^U,\bm\xi-\bm p^L)
\geq\bm0\}$.
\end{definition}


Figure~\ref{fig:dt} illustrates how a fixed tree structure, with optimizable thresholds $\bm p^U$ and $\bm p^L$, induces a parametric polyhedral partition of $\Xi$, where each leaf $k\in\mathcal K$ corresponds to a region $\Xi_k^{\bm\theta}$ defined by its root-to-leaf splitting rules.


\subsubsection{General Tree Partition Scheme}
The decision tree scheme restricts the target dimension to $\ell=2n$ and applies coordinate-wise splits. The following general tree scheme relaxes both restrictions at the cost of additional parameters.

\begin{definition}[General Tree Partition Scheme]
Given a tree with leaf set $\mathcal K$ and an even-dimensional target space $\mathbb R^\ell$, let $\bm s_k=(\bm s_k^U,\bm s_k^L)$ encode the root-to-leaf splits for each $k\in\mathcal K$. Define the target-space partition $\mathcal L=(\mathcal L_k)_{k\in\mathcal K}$ by $\mathcal L_k:=\{\bm x\in\mathbb R^\ell \mid\bm s_k\odot\bm x\geq\bm0\}$. Let $\bm\theta=(\bm P,\bm p^U,\bm p^L)$, where $\bm P\in\mathbb R^{(\ell/2)\times n}$ and $\bm p^U,\bm p^L\in\mathbb R^{\ell/2}$. The affine mapping is defined by $\bm T_{\bm\theta}\bm\xi+\bm\tau_{\bm\theta} :=(\bm P\bm\xi-\bm p^U,\bm P\bm\xi-\bm p^L)$, inducing the regions $\Xi_k^{\bm\theta}:= \{\bm\xi\in\Xi\mid \bm s_k\odot(\bm P\bm\xi-\bm p^U, \bm P\bm\xi-\bm p^L)\geq\bm0\}$.
\end{definition}

This scheme applies threshold splits to the transformed uncertainty $\bm P\bm\xi$ rather than directly to its coordinates, resembling optimal classification trees~\citep{bertsimas2017optimal}.
Although the orthant-based and decision-tree schemes are special cases of the general tree scheme, they offer distinct trade-offs in policy expressiveness, parameter dimensionality, and compatibility with the learning methods in Section~\ref{sec:learn}.

\subsection{Exact Solution Methods}
\label{sub-section-reformulation}
This subsection develops two exact reformulations of \eqref{eq:ppp_linear_kpart} under a given partition scheme: a master--subproblem method that iteratively generates worst-case scenarios, and a piecewise dualization method that exploits the polyhedral structure of $\Xi_k^{\bm\theta}$ to obtain a single-level formulation. 
Both reformulations serve as a basis for the scalable algorithm developed in Section~\ref{sec:learn}.



\subsubsection{Master--Subproblem Method}
\label{subsubsec:mastersub}
This method addresses the semi-infinite constraints in \eqref{eq:ppp_linear_kpart} by solving a master relaxation over a finite scenario set. Given a solution $(\bm y,\bm\theta)$, regional separation subproblems identify worst-case scenarios to be added to the master. The procedure terminates when the subproblems certify feasibility and optimality. The following theorem formalizes this approach.


\begin{restatable}{theorem}{mastersubproblem}
\label{thm:mastersub_reform}
Let $I$ be the index set of scenarios in  $\{\bm\xi_i\}_{i\in I}\subseteq \Xi$, and define binary variables $\lambda_{k,i}\in\{0,1\}$ to indicate whether $\bm\xi_i$ is assigned to region $\Xi_k^{\bm\theta}$. 
For sufficiently large constants $M$, $M_j$, and $M'$, Problem~\eqref{eq:ppp_linear_kpart} can be reformulated as the following master and subproblems:
\begin{subequations}
\label{eq:master}
  \vspace{-10pt}
\begin{align}
 z^{\mathrm{mst}} :=\quad 
\min_{\substack{\bm{y} \in \mathcal{Y}^{\mathcal{K}},z \in \mathbb{R} \\
 \bm \theta \in \Theta,\bm\lambda \in \{0,1\}^{\mathcal K \times I}}} &~ z\\
 \text{s.t.}\quad 
 &~ \iprod{\bm{Q}(\bm{y}_k), \bm{\xi}_i} + \bm{q}(\bm{y}_k) \le z + M\bigl(1-\lambda_{ki}\bigr), 
 \quad \forall i \in I,\ \forall k \in\mathcal{K},\label{eq:master-obj}\\
 &~ \iprod{\bm Q_j (\bm y_{k}), \bm\xi_i} + \bm q_j (\bm y_{k}) \le M_j\bigl(1-\lambda_{ki}\bigr), 
 \quad \forall j \in J,\ \forall k \in\mathcal{K},\ \forall i \in I,\label{eq:master-feas}\\
 &~\bm B_k \bigl(\bm T_{\bm \theta}\bm \xi_i + \bm \tau_{\bm\theta}\bigr) - \bm b_k 
 \le  M'\bigl(1-\lambda_{ki}\bigr)\boldsymbol{1}, 
 \quad \forall i \in I,\ \forall k \in\mathcal{K},\label{eq:master-in1}\\
  &~ \sum_{k\in\mathcal K} \lambda_{ki} = 1, \quad \forall i\in I.\label{eq:master-assign}
\end{align}
\end{subequations}
\par\noindent\begin{minipage}{\linewidth}
Given any solution $(\bm y, \bm\theta)$ from \eqref{eq:master}, define the following subproblems for every $k\in\mathcal K$ and $j\in J$:
\par\nobreak
\begin{subequations}
\label{eq:subproblems}
\par\nointerlineskip\vspace{4pt}
\noindent
\begin{minipage}[t]{0.48\textwidth}
\refstepcounter{equation}%
\label{eq:sub-obj}
\noindent\(\displaystyle
z^{\mathrm{opt}}_k := \
\begin{aligned}[t]
\max_{\bm\xi \in \Xi}&~\iprod{\bm{Q}(\bm{y}_k), \bm{\xi}} + \bm{q}(\bm{y}_k) \\
\text{s.t.}
&~\bm B_k \bigl(\bm T_{\bm \theta}\bm \xi + \bm \tau_{\bm\theta}\bigr)\le\bm b_k,
\end{aligned}
\)\hfill\textup{(\theequation)}
\end{minipage}%
\hfill
\begin{minipage}[t]{0.48\textwidth}
\refstepcounter{equation}%
\label{eq:sub-feas}
\noindent\(\displaystyle
z^{\mathrm{fea}}_{kj} := \
\begin{aligned}[t]
\max_{\bm\xi \in \Xi}&~\iprod{\bm Q_j (\bm y_{k}), \bm\xi} + \bm q_j (\bm y_{k}) \\
\text{s.t.}
&~\bm B_k \bigl(\bm T_{\bm \theta}\bm \xi + \bm \tau_{\bm\theta}\bigr)\le  \bm b_k,
\end{aligned}
\)\hfill\textup{(\theequation)}
\end{minipage}\par
\end{subequations}
\end{minipage}\par
\noindent The procedure halts when $z^{\mathrm{mst}} = \max_{k \in \mathcal K} z_k^{\mathrm{opt}}$ and $\max_{k \in \mathcal K, j \in J} z_{kj}^{\mathrm{fea}} \leq 0$ under some incumbent $(\bm y, \bm \theta)$.
\end{restatable}

Constraint~\eqref{eq:master-assign} assigns each scenario $\bm\xi_i$ to exactly one region, while \eqref{eq:master-in1} enforces membership $\bm\xi_i\in\Xi_k^{\bm\theta}$ whenever $\lambda_{ki}=1$. Constraints~\eqref{eq:master-obj} and \eqref{eq:master-feas} enforce the corresponding objective and feasibility requirements. The big-$M$ constants must be sufficiently large to deactivate these constraints when $\lambda_{ki}=0$. 
Although the separation subproblems can be solved in parallel, each new scenario introduces $K(1+|J|+\ell)$ constraints and $K$ binary variables in the master problem. This motivates decoupling partition learning and scenario assignments in Section~\ref{sec:learn}.

\subsubsection{Piecewise Dualization Method}
\label{subsubsec:dualization}

The explicit polyhedral description of each region $\Xi_k^{\bm\theta}$ enables us to dualize the regional worst-case objective and constraints. Applying LP duality region by region yields the following single-level reformulation of \eqref{eq:ppp_linear_kpart}.




\begin{restatable}{theorem}{dualization}
\label{thm:dualization}
Problem~\eqref{eq:ppp_linear_kpart} admits the following equivalent reformulation.
\begin{subequations}
\label{eq:dual-vi}
\begin{align}
    \min_{\substack{\bm y \in \mathcal Y^{\mathcal K},z\in\mathbb R\\ \bm \theta \in \Theta, \bm\pi \geq \bm 0, \bm\mu \geq \bm 0}}\quad
    & z\\
    \text{s.t.}\quad
    & \bm q(\bm y_k)+\iprod{\bm b,\bm\pi_k}+\iprod{\bm b_k-\bm B_k\bm\tau_{\bm\theta},\bm\mu_k}\le z, \quad \forall k\in\mathcal K,\\
    &\bm B^\top\bm\pi_k+\bm T_{\bm\theta}^\top\bm B_k^\top\bm\mu_k=\bm Q(\bm y_k),\quad \forall k\in\mathcal K,\\
    &\bm q_j(\bm y_k)+\iprod{\bm b,\bm\pi_{kj}}+\iprod{\bm b_k-\bm B_k\bm\tau_{\bm\theta},\bm\mu_{kj}} \le 0, \quad \forall k\in\mathcal K, j \in J,\\
    &\bm B^\top\bm\pi_{kj}+\bm T_{\bm\theta}^\top\bm B_k^\top\bm\mu_{kj}=\bm Q_j(\bm y_k), \quad \forall k\in\mathcal K, j \in J.
\end{align}
\end{subequations}
\end{restatable}

If an optimal solution exists, its partition is recovered directly from $\bm\theta^\star$. The reformulation is generally nonconvex due to bilinear terms involving the partition parameters $\bm\theta$ and dual variables $\bm\mu$, making direct solution by commercial solvers difficult to scale. Since these bilinear terms couple $\bm\theta$ with the remaining variables, a good initial partition $\bm\theta_0$ can be refined through alternating optimization over $(\bm y,z,\bm\pi,\bm\mu)$ and $\bm\theta$ \citep{grippo2000convergence}. We develop this approach further in the next section to obtain a scalable algorithm for polyhedral policies.


\section{A Scalable Learning-Based Framework for Finite Adaptability}
\label{sec:learn}


The exact methods in the previous section jointly optimize decisions and partition parameters, but become computationally expensive as problem size grows. To improve scalability, we propose the \emph{approximate-learn-parallel} (ALP) framework, which decouples partition learning from decision optimization. As illustrated in
Figure~\ref{fig:appx-learn-parallel}, ALP consists of four steps, detailed below.



\begin{figure}[!bp]
\FIGURE
{%
\centering
\begin{subfigure}[b]{\textwidth}
  \centering
\begingroup%
\newcommand{\alpFont}{\fontsize{8}{10}\selectfont\mdseries}%
\pgfmathsetlengthmacro{\alpUnit}{.9*\linewidth/15.65}%
\newcommand{\alpOutline}{%
  \draw[outline] (0,.55)--(.28,1.9)--(1.55,2.38)
    --(3.05,1.12)--(.88,0)--cycle;
  \node[anchor=north east,inner sep=0pt] at (2.92,.48) {$\Xi$};
}%
\newcommand{\alpCuts}{%
  \draw[cut] (.58,2.006)--(.88,0)
    (.76,2.081)--(1.72,.434)
    (.092,.994)--(2.22,1.817)
    (2.22,1.817)--(2.13,.645);
}%
\newcommand{\alpScenarios}{%
  \foreach \x/\y in {.54/1.71,.43/.94,.72/.48}
    \fill (\x,\y) circle[radius=1.7pt];
  \foreach \x/\y in {1.51/1.96,1.25/1.70,1.52/1.12}
    \path[draw=black,fill=black!65,line width=.4pt]
      (\x,\y) ++(0,2.5pt)--++(-2.2pt,-4.4pt)--++(4.4pt,0)--cycle;
  \foreach \x/\y in {.95/.93,1.17/.57}
    \path[draw=black,fill=black!65,line width=.4pt]
      (\x,\y) ++(-1.9pt,-1.9pt) rectangle ++(3.8pt,3.8pt);
  \foreach \x/\y in {1.81/1.21,2.30/1.42,1.78/.76}
    \path[draw=black,fill=black!65,line width=.4pt]
      (\x,\y) ++(0,2.6pt)--++(2.1pt,-2.6pt)
      --++(-2.1pt,-2.6pt)--++(-2.1pt,2.6pt)--cycle;
}%
\begin{tikzpicture}[x=\alpUnit,y=\alpUnit,
  every node/.style={font=\alpFont,text=black,inner sep=1pt},
  outline/.style={line width=.85pt,line join=round},
  cut/.style={line width=.65pt,black!80},
  flow/.style={draw=black,fill=white,line width=.7pt,line join=miter},
  heading/.style={align=left,anchor=north west,font=\alpFont},
  region/.style={inner sep=0pt}]
\begin{scope}[shift={(.183,0)},scale=.88]
  \alpOutline
  \alpScenarios
\end{scope}
\begin{scope}[shift={(4.383,0)},scale=.88]
  \alpOutline
  \alpCuts
  \alpScenarios
\end{scope}
\begin{scope}[shift={(8.583,0)},scale=.88]
  \alpOutline
  \alpCuts
  \node[region] at (.39,1.44) {$\Xi_1$};
  \node[region] at (.43,.65) {$\Xi_2$};
  \node[region] at (.84,1.61) {$\Xi_3$};
  \node[region] at (1.49,1.95) {$\Xi_4$};
  \node[region] at (1.72,1.21) {$\Xi_5$};
  \node[region] at (2.52,1.15) {$\Xi_6$};
  \node[region] at (1.13,.68) {$\Xi_7$};
\end{scope}
\begin{scope}[shift={(12.783,0)},scale=.88]
  \alpOutline
  \draw[cut] (.40,1.945)--(1.55,.346)
    (1.55,2.38)--(2.13,.645)
    (.092,.994)--(2.71,1.406);
  \node[region] at (.45,1.43) {$\Xi_1$};
  \node[region] at (.67,.60) {$\Xi_2$};
  \node[region] at (1.55,.88) {$\Xi_3$};
  \node[region] at (1.16,1.85) {$\Xi_4$};
  \node[region] at (2.06,1.74) {$\Xi_5$};
  \node[region] at (2.33,1.06) {$\Xi_6$};
\end{scope}
\foreach \x in {3.20,7.40,11.60}
  \path[flow] (\x,.95)--++(.51,0)--++(0,-.11)
    --++(.29,.22)--++(-.29,.22)--++(0,-.11)--++(-.51,0)--cycle;
\node[heading] at (.05,-.30)
  {Discrete approximation\\[3pt]$\triangleright$\enspace Exact method\\$\triangleright$\enspace Heuristic method};
\node[heading] at (4.25,-.30)
  {Partition learning\\[3pt]$\triangleright$\enspace SVMs\\$\triangleright$\enspace Decision trees};
\node[heading] at (8.45,-.30)
  {Parallel robust optimization};
\node[heading] at (12.65,-.30)
  {Partition enhancement\\[3pt](optional step)};
\path ([yshift=4pt]current bounding box.north west)
  rectangle (current bounding box.south east);
\end{tikzpicture}%
\endgroup%
%
\end{subfigure}
}
{\raggedright Four-step scalable approximate-learn-parallel (ALP) framework for finite adaptability.\label{fig:appx-learn-parallel}}
{First, a finite scenario set is initialized and labeled using either the exact or heuristic method. Second, a classifier, such as an SVM or a decision tree, is trained to learn a partition rule from these labels. Third, RO subproblems are solved independently over the resulting regions. Finally, the partition can be further enhanced using the dualization reformulation \eqref{eq:dual-vi}.}
\end{figure}

\subsection{Discrete Approximation}
\label{subsec:discrete_approx}
This step constructs a discrete partition to
support subsequent learning. Its quality depends on
representative scenarios $\{\bm\xi_i\}_{i\in I}$
and assignments that reflect the underlying objective
and constraint structure. Scenario selection depends
on the problem class. For example, under budget
uncertainty
$\Xi=\{\bm\xi\in[0,1]^n\mid
\iprod{\bm 1,\bm\xi}\le B\}$,
$B\in(0,n)$, a convex objective in $\bm\xi$
attains its worst case at an extreme point of $\Xi$.
If it is also componentwise nondecreasing,
the search can be restricted to extreme points
on $\iprod{\bm 1,\bm\xi}=B$.
For SP, scenarios are sampled
from the underlying distribution.

Given the scenarios, the discrete partition problem
is obtained by dropping \eqref{eq:master-in1}
from \eqref{eq:master}. When
$(\bm Q,\bm Q_j,\bm q,\bm q_j)$
have simple expressions in $\bm y_k$,
it can be solved using a standard MIP solver.
Since this step only approximates the partition,
solving it to optimality is unnecessary.
We therefore introduce a heuristic method in addition to the MIP formulation.

\subsubsection{Alternating Assignment--Optimization (AAO)}

Let
$f(\bm y_k,\bm\xi_i):=
\iprod{\bm Q(\bm y_k),\bm\xi_i}+\bm q(\bm y_k)$
and
$g_j(\bm y_k,\bm\xi_i):=
\iprod{\bm Q_j(\bm y_k),\bm\xi_i}+\bm q_j(\bm y_k)$.
Define $\delta_{\leq 0}(t)$ as $0$ if $t\leq 0$
and $+\infty$ otherwise.
Dropping \eqref{eq:master-in1} from
\eqref{eq:master} yields the following
assignment formulation:
\begin{subequations}
\label{eq:assign}
\begin{align}
  \min_{\bm\lambda\in\{0,1\}^{\mathcal K\times I}}
  &~\min_{\bm y\in\mathcal Y^{\mathcal K}}
  \max_{\substack{i\in I\\k\in\mathcal K}}
  h(\bm y_k,\bm\xi_i;\lambda_{ki})\\
  \text{s.t.}
  &~\sum_{k\in\mathcal K}\lambda_{ki}=1,
  \quad \forall i\in I,
\end{align}
\end{subequations}
where
\(
h(\bm y_k,\bm\xi_i;\lambda_{ki})
:=f(\bm y_k,\bm\xi_i)-M(1-\lambda_{ki})
+\sum_{j\in J}
\delta_{\leq 0}\bigl(
\lambda_{ki}g_j(\bm y_k,\bm\xi_i)
\bigr).
\)
We solve this approximately by alternating between scenario assignment
and $K$ independent optimization subproblems:
\begin{itemize}
  \item \emph{Initialization:}
  Obtain $K$ initial decisions $\{\bm y_k\}_{k\in\mathcal K}$
  heuristically, e.g., through solution elimination cuts,
  scenario clustering, or problem-specific methods
  such as $K$-shortest paths.

  \item \emph{Assignment (parallel):}
  For each $\bm\xi_i$, evaluate
  $\{h(\bm y_k,\bm\xi_i;1)\}_{k\in\mathcal K}$
  and assign it to the smallest minimizing index $k$,
  even if all policies are infeasible.

  \item \emph{Optimization (parallel):}
  Given $I_k:=\{i\in I:\lambda_{ki}=1\}$,
  solve independently for each nonempty $I_k$ the scenario-based robust problem
  \(
  \min_{\bm y_k\in\mathcal Y}
  \max_{i\in I_k}h(\bm y_k,\bm\xi_i;1).
  \)
\end{itemize}
The two steps repeat until the objective value
no longer decreases.
The optimization step is optimal for fixed
$\bm\lambda$. 
The following proposition establishes the optimality
of the assignment step and finite convergence of AAO.

\begin{restatable}{proposition}{assign}\label{prop:assig_optimal}
For any fixed policy $\{\bm y_k\}_{k\in\mathcal K}$, the assignment step is optimal for \eqref{eq:assign}. Moreover, AAO generates nonincreasing objective values and terminates in finitely many iterations at a blockwise optimal solution.
\end{restatable}

\subsection{Partition Learning}
\label{subsec:learn_partition}

The preceding discrete approximation step produces labeled scenarios $\{(\bm\xi_i,k_i)\}_{i\in I}$, which we use to learn a polyhedral partition of $\Xi$. We consider two classifiers corresponding to the partition schemes in Section~\ref{sec:model}: a linear multi-class SVM (one-vs-rest), whose separating hyperplanes define an orthant-based partition (Definition~\ref{defi:orthant}), and a decision tree, whose branching rules define a tree-based partition (Definition~\ref{defi:tree}). Other classifiers inducing polyhedral regions, such as XGBoost, can also be accommodated.
Unlike conventional predictive applications, the trained classifier is used to extract a polyhedral partition rather than to make predictions. Robustness is then enforced by solving the region-wise RO subproblems.

\subsection{Parallel Robust Optimization}
\label{subsec:parallel_opt}
Given the learned partition $\{\Xi_k^{\bm\theta}\}_{k\in\mathcal K'}$, we compute the recourse decisions by solving independent RO problems in parallel. For each $k\in\mathcal K'$, the subproblem is
\begin{subequations}\label{eq:R-Part}
\begin{align}
  \min_{\bm y\in\mathcal Y}
  \max_{\bm\xi\in\Xi_k^{\bm\theta}}
  &~\iprod{\bm Q(\bm y),\bm\xi}+\bm q(\bm y)\\
  \text{s.t.}\quad
  &\iprod{\bm Q_j(\bm y),\bm\xi}+\bm q_j(\bm y)
  \leq 0,
  \quad \forall\bm\xi\in\Xi_k^{\bm\theta},\ j\in J.
\end{align}
\end{subequations}
With $\bm\theta$ fixed, each subproblem admits the dual reformulation in \eqref{eq:dual-vi}. The resulting optimal decisions $\{\bm y_k^\star\}_{k\in\mathcal K'}$ define a polyhedral policy that selects $\bm y_k^\star$ for $\bm\xi\in\Xi_k^{\bm\theta}$.

 
\subsection{Partition Enhancement}
\label{subsec:partition_enhancement}
Starting from the partition learned in Step~2 and the decisions obtained in Step~3, we optionally refine the policy through alternating optimization of \eqref{eq:dual-vi}. Specifically, we alternate between optimizing $(\bm y,z,\bm\pi,\bm\mu)$ and $\bm\theta$ with the other set of variables fixed. Since all bilinear terms couple these two blocks, each update eliminates the bilinearity, allowing efficient solution by standard solvers. This procedure yields nonincreasing objective values and can further improve the learned partition.
As with other partition-based methods \citep{bertsimas2016adaptivepartitions,postek2016iterative}, ALP has no approximation guarantee for finite $K$. Nevertheless, Theorem~\ref{thm:convergence} establishes that polyhedral policies can approach the fully adjustable benchmark arbitrarily closely as $K$ increases.

\section{Extension to SP}
\label{sec:extsp}

Since polyhedral partitions are defined independently of the risk measure $\mathcal R_\Xi$, the framework naturally extends to SP.
Under a given partition scheme
$(\mathcal L,\bm T,\bm\tau)$,
let $h(\bm y,\bm\xi)$ denote the scenario cost,
equal to $f(\bm y,\bm\xi)$ when all constraints
$g_j(\bm y,\bm\xi)\leq 0$, $j\in J$, are satisfied,
and $+\infty$ otherwise.
The corresponding $K$-adaptable stochastic program is
\[
\min_{\substack{
\bm y:=(\bm y_k)_{k\in\mathcal K}
\in\mathcal Y^{\mathcal K}\\
\bm\theta\in\Theta}}
\sum_{k\in\mathcal K}
q_k^{\bm\theta}
\mathbb E\left[
h(\bm y_k,\bm\xi)
\mid \bm\xi\in\Xi_k^{\bm\theta}
\right],
\]
where $q_k^{\bm\theta}$ denotes the probability
of region $\Xi_k^{\bm\theta}$, with zero-probability
regions contributing zero.
Given i.i.d.\ samples $\{\bm\xi_i\}_{i\in I}$,
sample average approximation (SAA) yields
\begin{subequations}
\label{eq:saa}
\begin{align}
\min_{\substack{
\bm y\in\mathcal Y^{\mathcal K},
\bm z\in\mathbb R^{|I|}\\
\bm\theta\in\Theta,
\bm\lambda\in\{0,1\}^{\mathcal K\times I}}}
&~\frac{1}{|I|}\sum_{i\in I}z_i\\
\text{s.t.}\quad
& f(\bm y_k,\bm\xi_i)
\leq z_i+M(1-\lambda_{ki}),
\quad \forall i\in I,\ k\in\mathcal K,
\label{eq:saa01}\\
& g_j(\bm y_k,\bm\xi_i)
\leq M_j(1-\lambda_{ki}),
\quad \forall i\in I,\ k\in\mathcal K,\ j\in J,
\label{eq:saa-feas}\\
& \bm B_k
\bigl(\bm T_{\bm\theta}\bm\xi_i
+\bm\tau_{\bm\theta}\bigr)-\bm b_k
\leq M'(1-\lambda_{ki})\boldsymbol{1},
\quad \forall i\in I,\ k\in\mathcal K,
\label{eq:saa02}\\
& \sum_{k\in\mathcal K}\lambda_{ki}=1,
\quad \forall i\in I.
\label{eq:saa03}
\end{align}
\end{subequations}
The ALP framework can also be applied to
\eqref{eq:saa} to improve scalability.

\subsection{ALP Framework for Adaptable SP}
We extend the ALP framework in Section~\ref{sec:learn} to SP. In Step~1, given i.i.d.\ samples $\{\bm\xi_i\}_{i\in I}$, we solve \eqref{eq:saa} without \eqref{eq:saa02}, either exactly or through AAO with the same initialization and the following modified steps:
\begin{itemize}
    \item \emph{Assignment (parallel):}
    For each $\bm\xi_i$, evaluate
    $\{h(\bm y_k,\bm\xi_i)\}_{k\in\mathcal K}$
    and assign it to the first $k$ attaining the smallest value.

    \item \emph{Optimization (parallel):}
    Given $I_k:=\{i\in I:\lambda_{ki}=1\}$,
    solve independently for each nonempty $I_k$ the scenario-based stochastic problem
    \(
    \min_{\bm y_k\in\mathcal Y}
    \sum_{i\in I_k}h(\bm y_k,\bm\xi_i).
    \)
\end{itemize}
As in the robust setting, each step is optimal
with the other block fixed.
Upon termination, the resulting labels
$\{(\bm\xi_i,k_i)\}_{i\in I}$ are used in Step~2
to learn the partition parameter $\bm\theta$.
In Step~3, fixing $\bm\theta$, we solve
the regional stochastic programs
\(
\min_{\bm y\in\mathcal Y}
\mathbb E\left[
h(\bm y,\bm\xi)
\mid \bm\xi\in\Xi_k^{\bm\theta}
\right]
\)
for each $k\in\mathcal K'$ in parallel to obtain the final polyhedral policy.


\begin{remark}
Under complete recourse, \eqref{eq:saa} provides a sample average approximation of the stochastic program. Otherwise, feasibility over sampled scenarios does not guarantee feasibility over $\Xi$. The explicit polyhedral partition learned in Step~2 allows Step~3 to approximate the expected objective using i.i.d.\ samples while enforcing robust constraints over each region $\Xi_k^{\bm\theta}$. This yields robustly feasible policies with optimized sample-based expected performance.
\end{remark}

\section{Numerical Experiments}\label{sec:experiments}
This section reports the computational performance of the ALP framework on three classical testbeds for robust finite adaptability from the literature~\citep{hanasusanto2015k,bertsimas2016adaptivepartitions,subramanyam2020kadaptability}: (i) the shortest path problem with uncertain arc lengths (objective uncertainty), (ii) the capital budgeting problem with uncertain costs and profits (objective and single-constraint uncertainty), and (iii) the project management problem with uncertain task durations (size-dependent constraint uncertainty).
Section~\ref{sec:Ex-stochastic programing} then applies the stochastic extension of the ALP framework to the capital budgeting problem.

All experiments were conducted on a 2021 MacBook Pro equipped with an Apple M1 Max chip and 64\,GB of memory, running macOS Sonoma (version 14.6.1). Implementations were written in Python~3.13, with mixed-integer programs solved using Gurobi~13.0. The constraint feasibility tolerance was set to $\epsilon = 10^{-3}$ for accepting incumbent solutions. All other solver parameters were kept at their default values. We evaluate four variants of the ALP framework, distinguished by the algorithms employed in the first two steps:
    SVM--E: exact MIP in Step~1, SVM classifier in Step~2;
    DT--E: exact MIP in Step~1, decision tree classifier in Step~2;
    SVM--H: AAO heuristic in Step~1, SVM classifier in Step~2;
    DT--H: AAO heuristic in Step~1, decision tree classifier in Step~2.
For all these variants, Step~3 solves the region-wise robust (or stochastic) subproblems over the learned partition $\{\Xi_k^{\bm\theta}\}_{k\in\mathcal K'}$ in parallel, yielding solutions $\{\bm y_k\}_{k\in\mathcal K'}$.

For each instance, we test multiple values of $K \in \{2, 3, 4, 6, 8, 10\}$ and analyze the associated performance metrics. Specifically, we record the performance improvement~(\%), defined explicitly for each problem class, along with both the total parallel and sequential runtimes. In the parallel setting, the runtime of each parallel step (i.e., Step~1 of SVM--H and DT--H, Step~3 in all variants) is measured as the maximum completion time across all subroutines; in the sequential setting, it is measured as their sum. Both include the time spent on all three steps of ALP. A time limit of $1200$~seconds is imposed on Step 1 of all ALP variants. For the robust $K$-adaptability problems, we compare these variants against the baseline heuristic developed in~\cite{subramanyam2020kadaptability}, which, to the best of our knowledge, is among the most scalable algorithms currently available for general robust $K$-adaptability problems. We applied their publicly available code to our generated instances with a time limit of $1800$~seconds per value of~$K$. Since this baseline algorithm solves a sequence of RO problems with $K$ incremented by one at each iteration, its total runtime can exceed the per-$K$ limit of $1800$~seconds.

We also evaluate the partition enhancement step
(Step~4 of ALP), reporting its additional improvement
and runtime, with a five-minute time limit for each
alternating optimization subproblem.
Appendix~\ref{app:app_additional_results} provides
additional results for all experiments, comparing
parallel and sequential runtimes and examining
sensitivity to $K$ and the number of scenarios
used in Step~1.


\subsection{Shortest Paths}\label{shortes path}
This subsection studies the shortest path problem under the same setup as in~\citet{subramanyam2020kadaptability}. We analyze the results from two aspects: (i) scalability, assessed through performance improvement and runtime compared with the baseline algorithm, and (ii) additional improvement enabled by the partition enhancement step (Step~4 of ALP).


Let $G = (V, A)$ be a directed graph with node set~$V$, arc set~$A$, a source $s \in V$, and a terminal $t \in V$. The nominal arc lengths are denoted by $\bm{d}_0 \in \mathbb{R}_+^{|A|}$, and the realized lengths under a deviation scenario~$\bm{\xi}$ are given by
$
\bm{d}(\bm{\xi}) = (1 + \tfrac{1}{2}\,\bm{\xi}) \odot \bm{d}_0,
$
where $\bm{\xi}$ belongs to the budgeted uncertainty set
\(
\Xi = \{\, \bm{\xi} \in [0,1]^{|A|} \;|\; \langle \bm{1}, \bm{\xi} \rangle \le B \}
\)
for a fixed budget $B > 0$. The decision maker selects a menu of candidate $s$--$t$ paths \emph{here-and-now}, each represented by an incidence vector from the solution space
\[
\mathcal{Y} =
\left\{
\bm{y} \in \{0,1\}^{|A|}
\;\middle|\;
\sum_{(i,j) \in A} y_{ij} - \sum_{(j,i) \in A} y_{ji}
\geq \mathbb{I}[i = s] - \mathbb{I}[i = t],\; \forall\, i \in V
\right\},
\]
where $\mathbb{I}[\cdot]$ returns $1$ if the condition is satisfied and $0$ otherwise. After the uncertain arc lengths are realized, the least costly path among the candidates is implemented. This gives rise to the classic $K$-adaptability formulation
$
\min_{\bm{y} \in \mathcal{Y}^{\mathcal K}}\max_{\bm{\xi} \in \Xi}\min_{k \in \mathcal{K}}\;
\langle \bm{d}(\bm{\xi}),\, \bm{y}_k \rangle.
$ 
Since the partition learned in Step~2 of ALP may produce more regions than the initial~$K$, the resulting policy can exploit finer uncertainty partitions than a standard $K$-adaptability solution.

We consider graph sizes $N \in \{20, 50, 60, \ldots, 100\}$ and generate $10$ random directed graphs for each~$N$ as follows. Node coordinates are sampled independently and uniformly from $[0,10]^2$, and the nominal weight of each arc $(i,j) \in A$ is set to the Euclidean distance between its endpoints. The source~$s$ and terminal~$t$ are chosen as the pair of nodes with the maximum Euclidean distance. Starting from the complete directed graph on~$V$, we construct the arc set~$A$ by removing the $\lfloor 0.7(N^2 - N) \rfloor$ arcs with the largest nominal weights while preserving $s$--$t$ connectivity. The uncertainty budget is fixed at $B = 3$ for all instances.

\subsubsection{Scalability Comparison.}
\begin{table}[!b]
\TABLE
{\raggedright Performance comparison of ALP variants and the baseline algorithm on the shortest path problem. \label{tab:results_summary_nodes_edges}}
{\small\setlength{\tabcolsep}{4pt}\renewcommand{\arraystretch}{1.2}\centering
\begin{adjustbox}{max width=\textwidth,center}
\begin{tabular}{cc*{5}{c}c*{5}{c}}
\hline
\multirow{2}{*}{Nodes (\#)} & \multirow{2}{*}{Edges (\#)}
                            & \multicolumn{5}{c}{Improvement (\%)} &
& \multicolumn{5}{c}{Runtime (seconds)} \\
\cline{3-7}\cline{9-13}
& & SVM--E & DT--E & SVM--H & DT--H & Baseline & &
    SVM--E & DT--E & SVM--H & DT--H & Baseline\\
\hline
20  & 114  & 22.98 & 10.13 & 28.87 & 17.82 & \textbf{100.00} & & 30.59 & 30.64 &  5.25 &  5.24 & \textbf{0.05}\\
50  & 735  & 20.17 & 10.65 & 19.99 &  9.96 & \textbf{99.15}  & & 157.03& 26.05 & 32.75 & 32.45 & \textbf{11.44}\\
60  & 1062 & 19.98 &  8.28 & 17.70 &  9.93 & \textbf{98.41}  & & 152.45& \textbf{35.75} & 97.80 & 98.15 & 180.14\\
70  & 1449 & 18.48 &  9.81 & 19.35 & 11.13 & \textbf{97.88}  & & \textbf{46.01} & 46.18 & 140.23& 134.03& 521.63\\
80  & 1896 & 15.91 &  5.95 & 19.01 &  9.42 & \textbf{97.40}  & & \textbf{61.28} & 61.30 & 69.69 & 74.21 & 1336.45\\
90  & 2403 & 12.65 &  4.20 & 17.02 &  5.07 & \textbf{96.93}  & & 227.56& 227.30& \textbf{221.15} & 222.23& 1693.11\\
100 & 2970 & 12.51 &  2.09 & 15.41 &  4.03 & \textbf{96.45}  & & 253.17& 173.17& \textbf{108.07} & 108.81& 2377.88\\
\hline
\end{tabular}
\end{adjustbox}} 
{\begingroup\fontsize{8}{11}\selectfont\emph{Notes.}
Improvement is computed as $(z_{\mathrm{rob}}^\ast - z^\ast)/(z_{\mathrm{rob}}^\ast - z_{\mathrm{adj}}^\ast) \times 100$, where $z^\ast$, $z_{\mathrm{rob}}^\ast$, and $z_{\mathrm{adj}}^\ast$ denote the optimal values of the tested algorithm, static RO, and fully adjustable optimization. Step~1 uses $1500$ scenarios for the exact MIP and $8000$ for the AAO heuristic. For ALP, parallelized steps are timed by the maximum completion time across parallel tasks; corresponding sequential runtimes are reported in Appendix~\ref{app:app_additional_results}.
\par\endgroup}
\end{table}

Table~\ref{tab:results_summary_nodes_edges} summarizes the computational performance of the ALP framework (from Step~1 to Step~3, without Step~4) on the shortest path instances. For each graph size, we run every variant with $K \in \{2,3,4,6,8,10\}$ on ten randomly generated graphs and report the results associated with the~$K$ that yields the best average objective value~$z^\ast$. The column \emph{Improvement} is computed as $(z_{\mathrm{rob}}^\ast - z^\ast)/(z_{\mathrm{rob}}^\ast - z_{\mathrm{adj}}^\ast) \times 100$, where $z_{\mathrm{rob}}^\ast$ and $z_{\mathrm{adj}}^\ast$ denote the optimal values of the static RO model and the fully adjustable policy (obtained by dualizing the inner shortest-path problem). An improvement of $100\%$ indicates that the fully adjustable value is attained. The column \emph{Runtime} reports the average total wall-clock time under full parallelization.
The corresponding sequential runtimes, which assume all steps are solved in sequence without parallelization, are reported in Table~\ref{tab:results_summary_nodes_edges_time} in Appendix~\ref{app:app_additional_results}. For most instances, these sequential runtimes remain within approximately $120\%$ of the parallel counterparts.



Across all instance sizes, the baseline algorithm achieves the highest improvements, exceeding $95\%$ in every case, likely owing to its iterative scheme being particularly well suited to problems with objective-only uncertainty. Among the ALP variants, SVM--H performs best, with improvements ranging from $15\%$ to $29\%$.


In terms of runtime, the ALP variants are considerably more efficient than the baseline, with the largest instances solved within $300$~seconds, whereas the baseline runtime grows more rapidly with instance size. Among the ALP variants, the heuristic versions (SVM--H and DT--H) are generally faster than their exact counterparts (SVM--E and DT--E), particularly on larger instances. These results indicate that for robust shortest path problems, ALP provides a computationally efficient alternative that scales to larger networks, though it attains lower improvement than the baseline.


\subsubsection{Effect of Partition Enhancement.}
\begin{table}[!b]
\TABLE
{\raggedright Effect of partition enhancement (Step~4) on the shortest path problem. \label{tab:results_summary_nodes_edges_improvement}}
{\small\setlength{\tabcolsep}{4pt}\renewcommand{\arraystretch}{1.2}\centering
\begin{adjustbox}{max width=\textwidth,center}
\centering
\begin{tabular}{cc*{4}{c}c*{4}{c}}
\hline
\multirow{2}{*}{Nodes (\#)} & \multirow{2}{*}{Edges (\#)}
                            & \multicolumn{4}{c}{Additional Improvement (\%)} &
& \multicolumn{4}{c}{Additional Time (seconds)} \\
\cline{3-6}\cline{8-11}
& & SVM--E & DT--E & SVM--H & DT--H & &
    SVM--E & DT--E & SVM--H & DT--H \\
\hline
20  & 114  & +28.57 & \textbf{+40.87} & +18.91 & +29.92 & & \textbf{+2.04} & +11.64 & +2.70 & +5.94 \\
50  & 735  & +25.10 & \textbf{+27.50} & +20.38 & +25.83 & & +54.18 & +66.36 & +44.80 & \textbf{+32.86 } \\
60  & 1062 & +16.78 & \textbf{+27.84} & +17.49 & +20.57 & & +106.85 & +97.93 & +156.26 & \textbf{+65.54} \\
70  & 1449 & +21.77 & \textbf{+27.55} & +12.01 & +22.31 & & +282.40 & +263.85 & +345.23 & \textbf{+157.70} \\
80  & 1896 & +21.08 & +22.80 & +18.28 & \textbf{+25.53} & & +544.79 & \textbf{+147.42} & +1312.33 & +215.77 \\
90  & 2403 & +22.76 & +24.80 & +18.19 & \textbf{+25.03} & & +1530.08 & +543.55 & +868.60 & \textbf{+352.83} \\
100 & 2970 & +12.06 & \textbf{+27.51} & +14.20 & +18.44 & & +3067.45 & +823.65 & +1273.90 & \textbf{+626.02} \\
\hline
\end{tabular}
\end{adjustbox}}
{\begingroup\fontsize{8}{11}\selectfont\emph{Notes.}
Additional improvement and runtime from applying alternating optimization to the bilinear program~\eqref{eq:dual-vi}, warm-started with the Step~2 partition, for each ALP variant. For example, at $(20,114)$, SVM--E attains a total improvement of $22.98+28.57=51.55$ in $30.59+2.04=32.63$ seconds.
\par\endgroup}
\end{table}
Table~\ref{tab:results_summary_nodes_edges_improvement} reports the effect of the partition enhancement step (Step~4), which warm-starts the alternating optimization of the bilinear program~\eqref{eq:dual-vi} from the partition parameter~$\bm\theta$ obtained in Step~2. For each instance size, we report the additional improvement in percentage points and the corresponding additional runtime in seconds. For example, at the configuration $(20,114)$, SVM-E achieves $22.98 + 28.57=51.55$ total performance improvement with total runtime of $30.59 + 2.04=32.63$ seconds.

Across all graph sizes, Step~4 consistently improves the policies produced by ALP, yielding additional improvements ranging from $12$ to $40$ percentage points for every variant. Among all variants, DT--E achieves the largest incremental gains, potentially because its partition after Step~2 leaves more room for refinement. In terms of runtime, DT-based partitions lead to a more efficient Step~4, as the corresponding bilinear formulation~\eqref{eq:dual-vi} involves fewer partition parameters and is therefore smaller than its SVM-based counterpart.
These results indicate that Step~4 can effectively refine the learned partition and improve solution quality at the cost of additional computational time.

\subsection{Capital Budgeting}\label{capital budgeting}

Following the problem setup as in \citet{subramanyam2020kadaptability}, this subsection examines the capital budgeting problem.
A company allocates a fixed budget $B$ among $N$ projects whose cost vector $\bm{c}(\bm{\xi})$ and profit vector $\bm{r}(\bm{\xi})$ depend on a latent factor $\bm{\xi} \in \Xi:=[-1,1]^4$ via
$\bm c(\bm\xi)=(\bm 1+\tfrac{1}{2} \bm \Phi \bm\xi)\odot \bm c_{0}$ and 
$\bm r(\bm\xi)=(\bm 1+\tfrac{1}{2} \bm\Psi \bm\xi) \odot \bm r_{0}$
for given nominal vectors $\bm c_0$ and $\bm r_0$, and factor-loading matrices $\bm \Phi$ and $\bm \Psi$.
A discount factor $\rho = 0.8$ captures the fraction of potential profit realized upon investment. 
The finite-adaptable policy then optimizes over \(K\) candidate investment plans, selecting the most suitable one after observing \(\bm{\xi}\), defined as
\(
\max_{\bm{y} \in \mathcal{Y}^\mathcal K}
\min_{\bm{\xi} \in \Xi}
\max_{k \in \mathcal K}
\{
\rho\langle \bm{r}(\bm{\xi}), \bm{y}_k \rangle \, |\,
\langle \bm{c}(\bm{\xi}), \bm{y}_k \rangle \leq B
\},
\)
where $ \mathcal{Y} = \{0, 1\}^N $ and each $\bm{y}_k \in \mathcal Y$ indicates a project selection.
For each problem size $N \in \{5,30,50,100,1000,3000,5000\}$, we generate $10$ random instances. Nominal costs $\bm{c}_0$ are drawn uniformly from $[0,10]^N$, with $\bm{r}_0 = \bm{c}_0/5$. Each row of $\bm{\Phi}$ and $\bm{\Psi}$ is sampled uniformly from the unit simplex in $\mathbb{R}^4$, and the budget is set to $B = \tfrac{1}{2}\langle \bm{c}_0, \boldsymbol{1}\rangle$.

\subsubsection{Scalability Comparison.}
\begin{table}[!b]
\TABLE
{\raggedright Performance comparison of ALP variants and the baseline algorithm on the capital budgeting problem. \label{tab:results_summary}}
{\small\setlength{\tabcolsep}{4pt}\renewcommand{\arraystretch}{1.2}\centering
\begin{adjustbox}{center}
\begin{tabular}{c*{5}{c}c*{5}{c}}
\hline
\multirow{2}{*}{Projects (\#)}
& \multicolumn{5}{c}{Improvement (\%)}
&
& \multicolumn{5}{c}{Runtime (seconds)} \\
\cline{2-6}\cline{8-12}
& SVM--E & DT--E & SVM--H & DT--H & Baseline & &
  SVM--E & DT--E & SVM--H & DT--H & Baseline \\
\hline
5    & 82.19 & 48.53 & 41.45 & 41.54 & \textbf{127.02} & &
       4.29 & 16.78 & 0.21 & \textbf{0.14} & 0.15 \\
30   & \textbf{116.47} & 60.41 & 83.52 & 59.87 & 102.00 & &
       1240.56 & 1203.27 & 0.55 & \textbf{0.38} & 1912.99 \\
50   & \textbf{118.30} & 58.99 & 85.06 & 63.68 & 82.46 & &
       1220.19 & 1206.09 & 20.96 & \textbf{0.58} & 8696.41 \\
100  & \textbf{118.16} & 64.77 & 83.27 & 62.59 & -- & &
       1204.71 & 1213.52 & 33.74 & \textbf{0.94} & -- \\
1000 & \textbf{118.38} & 63.81 & 76.43 & 46.07 & -- & &
       1260.47 & 1253.39 & 5.97 & \textbf{5.88} & -- \\
3000 & \textbf{119.17} & 59.72 & 65.12 & 48.84 & -- & &
       1304.91 & 1360.52 & 18.69 & \textbf{18.07} & -- \\
5000 & \textbf{117.49} & 60.91 & 74.19 & 49.40 & -- & &
       1387.28 & 1471.68 & 30.12 & \textbf{30.07} & -- \\
\hline
\end{tabular}
\end{adjustbox}} 
{\begingroup\fontsize{8}{11}\selectfont\emph{Notes.}
Improvement is computed as $(z^\ast - z_{\mathrm{rob}}^\ast)/z_{\mathrm{rob}}^\ast \times 100$, where $z^\ast$ and $z_{\mathrm{rob}}^\ast$ denote the optimal values of the tested algorithm and static RO, respectively. Step~1 uses $500$ scenarios for both the exact MIP and the AAO heuristic. For ALP, parallelized steps are timed by the maximum completion time across parallel tasks; sequential runtimes are reported in Appendix~\ref{app:app_additional_results}.
\par\endgroup}
\end{table}
Table~\ref{tab:results_summary} summarizes the computational performance of the proposed ALP framework (from Step~1 to Step~3, without Step~4) on the capital budgeting instances. For each problem size $N$, we run every variant with $K \in \{2,3,4,6,8,10\}$ on ten randomly generated instances and report the results associated with the~$K$ that yields the best average objective value~$z^\ast$. 
Improvement is computed as $(z^\ast - z_{\mathrm{rob}}^\ast)/z_{\mathrm{rob}}^\ast \times 100$, differing from Table~\ref{tab:results_summary_nodes_edges} because the fully adjustable problem becomes difficult to solve when $\mathcal{Y}$ is discrete. 
Improvement can exceed 100\% because it is measured relative to the static RO objective value.
Runtime reports the average total wall-clock time under full parallelization. The corresponding sequential runtimes, reported in Table~\ref{tab:results_summary_seq} in Appendix~\ref{app:app_additional_results}, remain within roughly $200\%$ of their parallel counterparts across most instances.

Unlike the shortest path case, ALP attains both higher improvement and faster computation across all instances except the smallest ($N=5$). Even at $N=50$, the baseline achieves $82.46\%$ improvement in over two hours, whereas SVM--H reaches $85.06\%$ in under $25$ seconds. 
In terms of improvement, SVM-based methods outperform their DT counterparts, and exact MIP variants outperform their AAO heuristic counterparts. The heuristic variants, however, exhibit superior scalability, solving all instances in under one minute. Among all four ALP variants, SVM--H strikes the best balance between solution quality and computational efficiency.

\subsubsection{Effect of Partition Enhancement.}

Table~\ref{tab:results_summary_projects} reports the effect of the optional partition enhancement step (Step~4) in terms of additional improvement and the corresponding additional runtime. This partition enhancement step yields positive improvements across most tested instances, with the same trend of DT-based variants obtaining notably larger gains than their SVM counterparts. All instances are solved within minutes, and DT--H exhibits better efficiency on the larger instances. 
Overall, Step~4 provides a practical post-processing tool for improving finite-adaptable policies at modest additional computational cost, particularly when paired with DT-based variants.

\begin{table}[!tbp]
\TABLE
{\raggedright Effect of partition enhancement (Step~4) on the capital budgeting problem.\label{tab:results_summary_projects}}
{\small\setlength{\tabcolsep}{4pt}\renewcommand{\arraystretch}{1.2}\centering
\begin{adjustbox}{center}
\begin{tabular}{c*{4}{c}c*{4}{c}}
\hline
\multirow{2}{*}{Projects (\#)}
& \multicolumn{4}{c}{Additional Improvement (\%)}
&
& \multicolumn{4}{c}{Additional Time (seconds)} \\
\cline{2-5}\cline{7-10}
& SVM--E & DT--E & SVM--H & DT--H & &
  SVM--E & DT--E & SVM--H & DT--H \\
\hline
5    & +0.91 & \textbf{+34.90} & +0.00 & +32.00 & &
       \textbf{+0.16} & +0.22 & +8.37 & +0.31 \\
30   & +1.00 & \textbf{+9.45} & +0.32 & +7.41 & &
       +2.51 & \textbf{+0.96} & +1.51 & +1.36 \\
50   & +0.93 & \textbf{+7.23} & +0.22 & +3.36 & &
       +15.72 & \textbf{+0.56} & +2.24 & +4.17 \\
100  & +0.45 & +8.16 & +0.01 & \textbf{+8.64} & &
       +36.95 & +12.40 & +433.41 & \textbf{+10.59} \\
1000 & +0.01 & +2.42 & +2.03 & \textbf{+14.10} & &
       +21.59 & +16.64 & +634.60 & \textbf{+7.39} \\
3000 & +0.15 & +6.07 & +0.00 & \textbf{+11.37} & &
       +62.07 & +36.75 & +547.46 & \textbf{+ 26.03} \\
5000 & +0.03 & +20.09 & +0.00 & \textbf{+23.07} & &
       +91.48 & +96.78 & +787.11 & \textbf{+ 50.10} \\
\hline
\end{tabular}
\end{adjustbox}}
{\begingroup\fontsize{8}{11}\selectfont\emph{Notes.}
Additional improvement and runtime obtained by applying alternating optimization to the bilinear program~\eqref{eq:dual-vi}, warm-started with the Step~2 partition, for each ALP variant.
\par\endgroup}
\end{table}

\subsection{Project Management}\label{project management}

This subsection evaluates the algorithms on the project management problem~\citep{subramanyam2020kadaptability,wiesemann2012}. Let $G=(V,A)$ be a directed acyclic graph parameterized by $m \in \mathbb{Z}_+$, with node set $V:=[3m+1]$ and arc set $A:=\{(3l+1,3l+p),(3l+p,3l+4) \mid l=0,\ldots,m-1,\ p=2,3\}$. Nodes represent project tasks and arcs represent precedence relations. Each task $i \in V$ has uncertain duration $d_i(\bm{\xi})$ given by $d_{3l+2}(\bm{\xi})=\xi_{l+1}$, $d_{3l+3}(\bm{\xi})=1-\xi_{l+1}$ for $l=0,\ldots,m-1$, and $d_{3l+1}(\bm{\xi})=0$ for $l=0,\ldots,m$, with uncertainty set $\Xi=\{\bm{\xi}\in\mathbb{R}_+^m \mid \|\bm{\xi}-\boldsymbol{1}/2\|_1 \le 1/2\}$. The sink node $N:=3m+1$ represents project completion. A schedule $\bm{y}\in\mathcal{Y}=\mathbb{R}_+^N$ assigns start time $y_i$ to each task $i$. The $K$-adaptability problem minimizes the worst-case project makespan:
\[
\min_{\bm{y} \in \mathcal{Y}^\mathcal K}
\max_{\bm{\xi} \in \Xi}
\min_{k \in \mathcal K}
\bigl\{({\bm y}_k)_N\, \big|\,
({\bm y}_k)_j - ({\bm y}_k)_i \geq d_i(\bm \xi),\ \forall (i,j) \in A
\bigr\},
\]
We consider network sizes with \(m \in \{3,8,10,20,30,40,50\}\). As shown in \citet[Example~2.2]{wiesemann2012}, the optimal value of the static robust problem is \(m\). Unlike the previous problems, the number of constraints with right-hand-side uncertainty grows with the problem size.


\subsubsection{Scalability Comparison.}

Table~\ref{tab:results_summary_project_management} reports the scalability results. The corresponding sequential runtimes are reported in Table~\ref{tab:results_summary_project_management_time} in Appendix~\ref{app:app_additional_results}.
The baseline yields no improvement because, under constraint-only uncertainty, the initial robust policy has the same objective value across all scenarios. Improving the worst-case objective therefore requires better policies that collectively cover the entire uncertainty set, whereas the sequential heuristic adds one solution at a time while keeping the existing ones fixed. In contrast, all ALP variants achieve notable improvements in policy performance. Among the four variants, DT-based methods perform best on small and medium instances, while SVM--E achieves the greatest improvement on large instances. In terms of runtime, the heuristic variants complete within seconds, while their exact counterparts often reach the time limit at Step~1 of ALP.

\begin{table}[!tbp]
\TABLE
{\raggedright Performance comparison of ALP variants and the baseline algorithm on the project management problem. \label{tab:results_summary_project_management}}
{\small\setlength{\tabcolsep}{4pt}\renewcommand{\arraystretch}{1.2}\centering
\begin{adjustbox}{max width=\textwidth,center}
\begin{tabular}{cc*{5}{c}c*{5}{c}}
\hline
\multirow{2}{*}{\(m\)}
& \multirow{2}{*}{Projects (\#)}
& \multicolumn{5}{c}{Improvement (\%)}
&
& \multicolumn{5}{c}{Runtime (seconds)} \\
\cline{3-7}\cline{9-13}
& & SVM--E & DT--E & SVM--H & DT--H & Baseline & &
  SVM--E & DT--E & SVM--H & DT--H & Baseline \\
\hline
3    & 10    &  0.00  &  \textbf{18.03} &  0.07   & 16.23 & 0.00 & &
      1211.21  & 1080.53  & 0.93 & \textbf{0.91} & 0.00 \\
8   & 25   &   1.21  &  8.77 & 8.19 & \textbf{11.93} & 0.00 & &
      1212.12   &  467.18  & 1.89 & \textbf{1.74} & 0.13 \\
10   & 31   & 0.81   &  7.02  &  7.14 & \textbf{9.77} & 0.00 & &
        1236.96   &  1207.55  &   2.28 & \textbf{2.07} & 0.12 \\
20  & 61  &   2.93  &  1.61  & 4.22   & \textbf{4.87} & 0.00 & &
      1273.80    &  1215.18  &  4.22  & \textbf{3.86}  & 0.32 \\
30 & 91 &  \textbf{3.25}   &  0.73  &  2.99 & 3.21 & 0.00 & &
     1310.20     &   1310.09  &   6.09  & \textbf{5.71} & 1.32 \\
40 & 121 &   \textbf{3.20}  &  0.82   &  2.31 & 2.35 & 0.00 & &
        1348.77  &  1318.73  &   8.11  & \textbf{7.07} & 1.67 \\
50 & 151 &   \textbf{3.68}   &  0.61  &  1.88 & 1.89  & 0.00 & &
      1312.53   &   1256.97  &  10.32   & \textbf{8.98} & 1.39 \\
\hline
\end{tabular}
\end{adjustbox}} 
{\begingroup\fontsize{8}{11}\selectfont\emph{Notes.}
Improvement is computed as $(z_{\mathrm{rob}}^\ast - z^\ast)/z_{\mathrm{rob}}^\ast \times 100$, where $z^\ast$ and $z_{\mathrm{rob}}^\ast$ denote the optimal values of the tested algorithm and static RO. Step~1 uses $2000$ scenarios for all variants. For ALP, parallelized steps are timed by the maximum completion time across parallel tasks; sequential runtimes are reported in Appendix~\ref{app:app_additional_results}.
Baseline runtimes are not bolded as no improvement is achieved.
\par\endgroup}
\end{table}

\subsubsection{Effect of Partition Enhancement.}

Table~\ref{tab:results_summary_project_management_bilinear} reports the effect of the optional partition refinement step (Step~4) on the project management instances. Step~4 provides additional improvement mainly for the DT-based variants, while the SVM-based variants show little or no gain. This difference arises because SVM-based partitions typically contain more regions, resulting in larger formulations that cannot be optimized within the five-minute time limit. In contrast, DT-based partitions keep the problem size smaller, allowing Step~4 to remain computationally effective, with DT--E obtaining the largest gains for most instance sizes. DT-based variants also require less additional runtime than their SVM counterparts, again due to the smaller problem sizes.

\begin{table}[!h]
\TABLE
{\raggedright Effect of partition enhancement (Step~4) on the project management problem.\label{tab:results_summary_project_management_bilinear}}
{\small\setlength{\tabcolsep}{4pt}\renewcommand{\arraystretch}{1.2}\centering
\begin{adjustbox}{center}
\begin{tabular}{cc*{4}{c}c*{4}{c}}
\hline
\multirow{2}{*}{$m$}
& \multirow{2}{*}{Projects (\#)}
& \multicolumn{4}{c}{Additional Improvement (\%)}
&
& \multicolumn{4}{c}{Additional Time (seconds)} \\
\cline{3-6}\cline{8-11}
& & SVM--E & DT--E & SVM--H & DT--H & &
  SVM--E & DT--E & SVM--H & DT--H \\
\hline
  3  &    10  & +0.32 & \textbf{+52.05} &   +0.00 & +25.69 & &
      \textbf{+0.53} & +0.80 & +8.37 & +0.90 \\
    8 &  25    & +0.22 & \textbf{+6.99} & +0.00 & +6.30 & &
      +3.07  & \textbf{+2.02} & +76.27 & +4.82 \\
  10   &  31    & +0.00 & \textbf{+8.53} & +0.00 & +2.16 & &
      +110.41  & \textbf{+5.18} & +91.14 & +7.15 \\
   20  &  61    & +0.00 & \textbf{+4.26} & +0.00 & +0.00 & &
      +283.54  & +22.10 & +319.70 & \textbf{+2.06} \\
  30   &   91   & +0.00 & \textbf{+4.25} & +0.00  & +1.15 & &
     +237.02   & +36.82 & +748.01 & \textbf{+7.44} \\
  40   &   121   & +0.00 & \textbf{+3.29} & +0.00  & +1.86 & &
     +73.12   & +48.96 & +296.60 & \textbf{+15.48} \\
   50  &   151   & +0.00 & +1.56 & +0.00  & \textbf{+2.32} & &
      +93.62  & +90.19 & +509.51 & \textbf{+17.33} \\
\hline
\end{tabular}
\end{adjustbox}}
{\begingroup\fontsize{8}{11}\selectfont\emph{Notes.}
Additional improvement and runtime obtained by applying alternating optimization to the bilinear program~\eqref{eq:dual-vi}, warm-started with the Step~2 partition, for each ALP variant.
\par\endgroup}
\end{table}

\subsection{ALP in SP}\label{sec:Ex-stochastic programing}
This subsection evaluates the stochastic version of ALP, introduced in Section~\ref{sec:extsp}, on the capital budgeting problem. Given their scalable performance in the robust setting, we focus on the SVM--H and DT--H variants. For each problem size $N \in \{5,30,50,100,1000,3000,5000\}$, five instances are generated following the setup in Subsection~\ref{capital budgeting}. We note that existing methods on stochastic $K$-adaptability \citep{buchheim2019k,malaguti2022k,donninik} pursue exact solutions for discrete uncertainty sets, while the proposed ALP framework targets approximate finite-adaptable policies for general uncertainty sets with an emphasis on scalability.


Table~\ref{tab:results_summary_stochastic_parallel} reports the computational results. For each problem size $N$, every variant is evaluated over $K \in \{2,3,4,6,8,10\}$, and we report results for the $K$ attaining the best average objective value $z^\ast$. Improvement is computed as $(z^\ast - z_{\mathrm{sto}}^\ast)/z_{\mathrm{sto}}^\ast \times 100$, where $z_{\mathrm{sto}}^\ast$ denotes the optimal value of the SP benchmark. Both parallel and sequential runtimes are presented.

\begin{table}[!b]
\TABLE
{\raggedright Performance comparison of the stochastic extension of ALP variants on the capital budgeting problem.\label{tab:results_summary_stochastic_parallel}}
{\small\setlength{\tabcolsep}{4pt}\renewcommand{\arraystretch}{1.2}\centering
\begin{adjustbox}{center}
\begin{tabular}{c*{2}{c}c*{2}{c}c*{2}{c}}
\hline
\multirow{2}{*}{Projects (\#)}
& \multicolumn{2}{c}{Improvement (\%)}
& & \multicolumn{2}{c}{Parallel Time (seconds)}
& & \multicolumn{2}{c}{Sequential Time (seconds)} \\
\cline{2-3}\cline{5-6}\cline{8-9}
& SVM--H & DT--H & & SVM--H & DT--H & & SVM--H & DT--H \\
\hline
5    & 22.53 & \textbf{24.23} & & \textbf{0.95}   & 1.60   & & \textbf{1.29}   & 9.14 \\
30   & 24.90 & \textbf{27.90} & & 5.06   & \textbf{4.82}   & & \textbf{14.30}  & 45.89 \\
50   & 25.12 & \textbf{28.87} & & 12.98  & \textbf{7.39}   & & 72.00  & \textbf{71.99} \\
100  & 25.18 & \textbf{29.35} & & 13.06  & \textbf{10.62}  & & \textbf{83.15}  & 118.07 \\
1000 & 24.84 & \textbf{29.15} & & 81.63  & \textbf{79.91}  & & \textbf{251.76} & 1183.31 \\
3000 & 22.29 & \textbf{28.28} & & \textbf{233.58} & 235.18 & & \textbf{832.39} & 3289.97 \\
5000 & 22.54 & \textbf{27.98} & & \textbf{385.23} & 398.66 & & \textbf{1273.95} & 5083.45 \\
\hline
\end{tabular}
\end{adjustbox}}
{\begingroup\fontsize{8}{11}\selectfont\emph{Notes.}
Improvement is computed as $(z^\ast - z_{\mathrm{sto}}^\ast)/z_{\mathrm{sto}}^\ast \times 100$, where $z^\ast$ and $z_{\mathrm{sto}}^\ast$ denote the optimal values of the tested algorithm and the SP benchmark, respectively. Step~1 uses $2000$ scenarios for the AAO heuristic. Parallel times report the maximum completion time across parallel tasks; sequential times report the sum of these task durations.
\par\endgroup}
\end{table}

Across all instance sizes, DT--H consistently achieves higher improvement than SVM--H, ranging from $24\%$ to $29\%$ compared with $22\%$ to $25\%$ for SVM--H. Moreover, improvement remains stable as $N$ increases, indicating that the stochastic ALP framework maintains solution quality even on the largest instances.
In terms of runtime, both variants remain computationally tractable as $N$ increases. Under parallelization, they perform similarly, with DT--H slightly faster on small-to-medium instances and SVM--H on the largest, and all runtimes remaining below $400$ seconds. In contrast, sequential runtimes grow noticeably with $N$, especially for DT--H, whose learning step refines the uncertainty set into a larger number of regions. Overall, the stochastic extension of ALP achieves stable improvement while remaining efficient, particularly under parallel execution.


\subsection{Summary of Numerical Insights}
Beyond computational scalability, our experiments offer further insights into finite adaptability. When uncertainty affects only the objective (e.g., the shortest path problem), exact worst-case evaluation largely determines solution quality, and algebraic approaches such as \cite{subramanyam2020kadaptability} hold a natural advantage. 
When uncertainty also affects constraints, as in the capital budgeting and project management problems, the proposed partition-based approach achieves stronger objective values and scales to substantially larger instances.
These observations are also consistent with existing results, such as \cite{hanasusanto2015k}, which shows that worst-case evaluation becomes intractable under constraint uncertainty, and \cite{han2023finiteadaptability}, which characterizes the non-convex structure of optimal partitions.

\section{Conclusion}\label{sec:conclusion}
This paper developed a scalable framework for finite adaptability in both robust and stochastic settings. We established that polyhedral finite-adaptable policies converge to the fully adjustable benchmark under mild regularity conditions, motivating a parametric polyhedral partition framework that offers flexible design of partition schemes. Building on this, we proposed the ALP framework, which decouples partition construction from decision optimization by learning the partition from a discrete approximation and solving the resulting region-wise subproblems in parallel. The learned partition can be further refined through alternating optimization on a piecewise-dualization reformulation. Computational experiments on shortest path, capital budgeting, and project management problems demonstrated that the proposed approach scales to instance sizes that are challenging for existing approaches, while producing solutions that are independently optimized within each region.
 
Several directions for future exploration emerge naturally from this work.
First, extending the computational framework to two-stage problems with explicit first-stage decisions and to multistage settings would broaden its applicability.
Second, the partition-based structure is not inherently tied to a specific risk measure, suggesting natural extensions to DRO. 
Third, investigating richer partition schemes, such as those induced by ensemble trees or piecewise-linear neural network classifiers, may improve partition quality while maintaining tractable polyhedral representations.


\clearpage
\bibliographystyle{plainnat}
\bibliography{bibi/references}

@article{julien2025machine,
  title={Machine learning for k-adaptability in two-stage robust optimization},
  author={Julien, Esther and Postek, Krzysztof and Birbil, {\c{S}} {\.I}lker},
  journal={INFORMS Journal on Computing},
  volume={37},
  number={3},
  pages={644--665},
  year={2025},
  publisher={INFORMS}
}

@misc{kedadsidhoum2025,
      title={Finite adaptability in two-stage robust optimization: asymptotic optimality and tractability}, 
      author={Safia Kedad-Sidhoum and Anton Medvedev and Frédéric Meunier},
      year={2025},
      eprint={2305.05399},
      archivePrefix={arXiv},
      primaryClass={math.OC},
}

@article{bertsimas2011geometric,
  title={A geometric characterization of the power of finite adaptability in multistage stochastic and adaptive optimization},
  author={Bertsimas, Dimitris and Goyal, Vineet and Sun, Xu Andy},
  journal={Mathematics of Operations Research},
  volume={36},
  number={1},
  pages={24--54},
  year={2011},
  publisher={INFORMS}
}

@article{ben2004adjustable,
  title={Adjustable robust solutions of uncertain linear programs},
  author={Ben-Tal, Aharon and Goryashko, Alexander and Guslitzer, Elana and Nemirovski, Arkadi},
  journal={Mathematical Programming},
  volume={99},
  number={2},
  pages={351--376},
  year={2004},
  publisher={Springer}
}

@article{chen2008linear,
  title={A linear decision-based approximation approach to stochastic programming},
  author={Chen, Xin and Sim, Melvyn and Sun, Peng and Zhang, Jiawei},
  journal={Operations Research},
  volume={56},
  number={2},
  pages={344--357},
  year={2008},
  publisher={INFORMS}
}

@article{yanikouglu2019survey,
  title={A survey of adjustable robust optimization},
author={Yanikoglu, Ihsan and Gorissen, Bram L and den Hertog, Dick},
  journal={European Journal of Operational Research},
  volume={277},
  number={3},
  pages={799--813},
  year={2019},
  publisher={Elsevier}
}

@article{malaguti2022k,
  title={K-adaptability in stochastic optimization},
  author={Malaguti, Enrico and Monaci, Michele and Pruente, Jonas},
  journal={Mathematical Programming},
  volume={196},
  number={1},
  pages={567--595},
  year={2022},
  publisher={Springer}
}

@article{buchheim2019k,
  title={K-adaptability in stochastic combinatorial optimization under objective uncertainty},
  author={Buchheim, Christoph and Pruente, Jonas},
  journal={European Journal of Operational Research},
  volume={277},
  number={3},
  pages={953--963},
  year={2019},
  publisher={Elsevier}
}

@book{ben2009robust,
  title     = {Robust Optimization},
  author    = {Ben-Tal, Aharon and El Ghaoui, Laurent and Nemirovski, Arkadi},
  year      = {2009},
  publisher = {Princeton University Press}
}

@article{bertsimas2010finite,
  author    = {Bertsimas, Dimitris and Caramanis, Constantine},
  title     = {Finite adaptability in multistage linear optimization},
  journal   = {IEEE Transactions on Automatic Control},
  volume    = {55},
  number    = {12},
  pages     = {2751--2766},
  year      = {2010}
}

@article{bertsimas2010optimality,
  author    = {Bertsimas, Dimitris and Iancu, Dan A. and Parrilo, Pablo A.},
  title     = {Optimality of affine policies in multistage robust optimization},
  journal   = {Mathematics of Operations Research},
  volume    = {35},
  number    = {2},
  pages     = {363--394},
  year      = {2010}
}

@article{bertsimas2017optimal,
  author    = {Bertsimas, Dimitris and Dunn, Jack},
  title     = {Optimal classification trees},
  journal   = {Machine Learning},
  volume    = {106},
  number    = {7},
  pages     = {1039--1082},
  year      = {2017},
  publisher = {Springer}
}

@article{grippo2000convergence,
  author    = {Grippo, Luigi and Sciandrone, Marco},
  title     = {On the convergence of the block nonlinear Gauss--Seidel method under convex constraints},
  journal   = {Operations Research Letters},
  volume    = {26},
  number    = {3},
  pages     = {127--136},
  year      = {2000},
  publisher = {Elsevier}
}

@article{han2023finiteadaptability,
  author    = {Han, Eojin and Bandi, Chaithanya and Nohadani, Omid},
  title     = {On finite adaptability in two-stage distributionally robust optimization},
  journal   = {Operations Research},
  volume    = {71},
  number    = {1},
  pages     = {88--106},
  year      = {2023}
}

@article{hanasusanto2015k,
  author    = {Hanasusanto, Grani A. and Kuhn, Daniel and Wiesemann, Wolfram},
  title     = {K-Adaptability in Two-Stage Robust Binary Programming},
  journal   = {Operations Research},
  volume    = {63},
  number    = {4},
  pages     = {877--891},
  year      = {2015}
}

@book{shapiro2021lectures,
  author    = {Shapiro, Alexander and Dentcheva, Darinka and Ruszczynski, Andrzej},
  title     = {Lectures on Stochastic Programming: Modeling and Theory},
  year      = {2021},
  publisher = {SIAM}
}

@article{subramanyam2020kadaptability,
  author    = {Subramanyam, Anirudh and Gounaris, Chrysanthos E. and Wiesemann, Wolfram},
  title     = {K-adaptability in two-stage mixed-integer robust optimization},
  journal   = {Mathematical Programming},
  volume    = {184},
  number    = {1},
  pages     = {177--219},
  year      = {2020}
}

@article{kuhn2011primal,
  author    = {Kuhn, Daniel and Wiesemann, Wolfram and Georghiou, Angelos},
  title     = {Primal and dual linear decision rules in stochastic and robust optimization},
  journal   = {Mathematical Programming},
  volume    = {130},
  number    = {1},
  pages     = {177--209},
  year      = {2011}
}

@article{georghiou2015generalized,
  author    = {Georghiou, Angelos and Wiesemann, Wolfram and Kuhn, Daniel},
  title     = {Generalized decision rule approximations for stochastic programming via liftings},
  journal   = {Mathematical Programming},
  volume    = {152},
  number    = {1},
  pages     = {301--338},
  year      = {2015}
}

@article{bertsimas2012power,
  author    = {Bertsimas, Dimitris and Goyal, Vineet},
  title     = {On the power and limitations of affine policies in two-stage adaptive optimization},
  journal   = {Mathematical Programming},
  volume    = {134},
  number    = {2},
  pages     = {491--531},
  year      = {2012}
}

@article{georghiou2020primal,
  author    = {Georghiou, Angelos and Tsoukalas, Angelos and Wiesemann, Wolfram},
  title     = {A primal--dual lifting scheme for two-stage robust optimization},
  journal   = {Operations Research},
  volume    = {68},
  number    = {2},
  pages     = {572--590},
  year      = {2020}
}

@article{ben2020tractable,
  author    = {Ben-Tal, Aharon and El Housni, Omar and Goyal, Vineet},
  title     = {A tractable approach for designing piecewise affine policies in two-stage adjustable robust optimization},
  journal   = {Mathematical Programming},
  volume    = {182},
  number    = {1},
  pages     = {57--102},
  year      = {2020}
}

@article{bertsimas2016adaptivepartitions,
  author    = {Bertsimas, Dimitris and Dunning, Iain},
  title     = {Multistage Robust Mixed-Integer Optimization with Adaptive Partitions},
  journal   = {Operations Research},
  volume    = {64},
  number    = {4},
  pages     = {980--998},
  year      = {2016}
}

@article{postek2016iterative,
  author    = {Postek, Krzysztof and den Hertog, Dick},
  title     = {Multistage Adjustable Robust Mixed-Integer Optimization via Iterative Splitting of the Uncertainty Set},
  journal   = {INFORMS Journal on Computing},
  volume    = {28},
  number    = {3},
  pages     = {553--574},
  year      = {2016}
}

@article{buchheim2017minmaxmin,
  author    = {Buchheim, Christoph and Kurtz, Jannis},
  title     = {Min--max--min robust combinatorial optimization},
  journal   = {Mathematical Programming},
  volume    = {163},
  number    = {1},
  pages     = {1--23},
  year      = {2017}
}

@article{buchheim2018complexity,
  author    = {Buchheim, Christoph and Kurtz, Jannis},
  title     = {Complexity of min--max--min robustness for combinatorial optimization under discrete uncertainty},
  journal   = {Discrete Optimization},
  volume    = {28},
  pages     = {1--15},
  year      = {2018}
}

@article{kurtz2024approximation,
  author    = {Kurtz, Jannis},
  title     = {Approximation Guarantees for Min-Max-Min Robust Optimization and $k$-Adaptability Under Objective Uncertainty},
  journal   = {SIAM Journal on Optimization},
  volume    = {34},
  number    = {2},
  pages     = {1538--1563},
  year      = {2024}
}

@article{goerigk2020minmaxmin,
  author    = {Goerigk, Marc and Kurtz, Jannis and Poss, Michael},
  title     = {Min--max--min robustness for combinatorial problems with discrete budgeted uncertainty},
  journal   = {Discrete Applied Mathematics},
  volume    = {285},
  pages     = {707--725},
  year      = {2020}
}

@article{bertsimas2015design,
  author    = {Bertsimas, Dimitris and Georghiou, Angelos},
  title     = {Design of Near Optimal Decision Rules in Multistage Adaptive Mixed-Integer Optimization},
  journal   = {Operations Research},
  volume    = {63},
  number    = {3},
  pages     = {610--627},
  year      = {2015}
}

@misc{donninik,
  author       = {Donnini, Federica and De Santis, Marianna and Kurtz, Jannis},
  title        = {K-adaptability for two-stage stochastic optimization},
  year         = {2025},
  howpublished = {Optimization Online}
}

@article{hanasusanto2016k,
  title={K-adaptability in two-stage distributionally robust binary programming},
  author={Hanasusanto, Grani A and Kuhn, Daniel and Wiesemann, Wolfram},
  journal={Operations Research Letters},
  volume={44},
  number={1},
  pages={6--11},
  year={2016},
  publisher={Elsevier}
}

@article{daryalal2024lagrangian,
  title={Lagrangian dual decision rules for multistage stochastic mixed-integer programming},
  author={Daryalal, Maryam and Bodur, Merve and Luedtke, James R},
  journal={Operations Research},
  volume={72},
  number={2},
  pages={717--737},
  year={2024},
  publisher={INFORMS}
}

@article{thoma2026note,
  title={A note on piecewise affine decision rules for robust, stochastic, and data-driven optimization},
  author={Thoma, Simon and Schiffer, Maximilian and Wiesemann, Wolfram},
  journal={Operations Research},
  year={2026},
  publisher={INFORMS}
}

@article{han2025nonlinear,
  author  = {Han, Eojin and Nohadani, Omid},
  title   = {Nonlinear Decision Rules Made Scalable by Nonparametric Liftings},
  journal = {Management Science},
  volume  = {71},
  number  = {4},
  pages   = {3449--3471},
  year    = {2025}
}

@article{wiesemann2012,
  author  = {Wiesemann, Wolfram and Kuhn, Daniel and Rustem, Ber\c{c}},
  title   = {Robust resource allocations in temporal networks},
  journal = {Mathematical Programming},
  volume  = {135},
  number  = {1},
  pages   = {437--471},
  year    = {2012}
}

\clearpage
\appendix

\section{Mathematical Proofs}


\noindent\textbf{Proof of Proposition~\ref{prop:piecewise}.}
Let $(\bm y,\mathcal P)$ be an optimal solution of~\eqref{eq:kpart}, where $\bm y=(\bm y_k)_{k\in\mathcal K}$. For fixed $\bm y$, assigning each $\bm\xi\in\Xi$ to a response minimizing $f(\bm y_k,\bm\xi)$ cannot increase the objective by monotonicity of $\mathcal R_\Xi$.
If $f(\bm y_k,\cdot)$ is affine for every
$k\in\mathcal K$, the regions
\[
\Xi_k=\{\bm\xi\in\Xi:
f(\bm y_k,\bm\xi)\leq
f(\bm y_{k'},\bm\xi),\ \forall k'\neq k\}
\]
are intersections of $\Xi$ with finitely many
half-spaces, and hence are polyhedral relative
to $\Xi$. Ties are resolved by the selection
rule $\sigma_{\mathcal P}(\cdot)$ without
affecting the objective.
If $f(\bm y_k,\cdot)$ is piecewise affine with
finitely many polyhedral pieces, take a common
polyhedral refinement of these pieces over all
$k\in\mathcal K$. On each refined cell,
the differences
$f(\bm y_k,\bm\xi)-f(\bm y_{k'},\bm\xi)$
are affine. Subdividing each cell along the
corresponding comparison hyperplanes yields
finitely many polyhedral cells on which a
minimizing response can be selected.
Grouping cells assigned to the same response
gives an optimal partition whose regions are
finite unions of polyhedra.
\pfend

The following two lemmas are needed for Theorem~\ref{thm:convergence}.

\begin{restatable}{lemma}{riskmeasure}
\label{lem:riskmeasure}
Every risk measure $\mathcal R_\Xi$ is 1-Lipschitz under the supremum norm, i.e.,
$|\mathcal R_\Xi(f)-\mathcal R_\Xi(g)|
\leq \|f-g\|_\infty$
for any real-valued functions $f,g$.
\end{restatable}

\noindent\textbf{Proof.}
Let $\delta:=\|f-g\|_\infty$. Then
$g-\delta\leq f\leq g+\delta$.
By monotonicity and translation equivariance,
\[
\mathcal R_\Xi(g)-\delta
=\mathcal R_\Xi(g-\delta)
\leq\mathcal R_\Xi(f)
\leq\mathcal R_\Xi(g+\delta)
=\mathcal R_\Xi(g)+\delta.
\]
Hence,
$|\mathcal R_\Xi(f)-\mathcal R_\Xi(g)|
\leq\|f-g\|_\infty$.
\pfend

\begin{restatable}{lemma}{polygon}
\label{lem:convergence}
Let $\Xi\subseteq\mathbb R^n$ be compact and contained in a hypercube of side length $D>0$, and consider
$\pi(\bm\xi)\in\arg\min_{\bm y\in\mathcal Y}
f(\bm y,\bm\xi)$, where
$f:\mathcal Y\times\Xi\to\mathbb R\cup\{+\infty\}$.
Suppose $\phi_\pi$ is finite on $\Xi$.
For every $\epsilon>0$, define the local robust value
$\psi_\epsilon(\bm\xi)
:=\inf_{\bm y\in\mathcal Y}
\sup_{\bm\xi'\in\bar B_{\epsilon,\bm\xi}}
f(\bm y,\bm\xi')$,
where $\bar B_{\epsilon,\bm\xi}
:=\{\bm\xi'\in\Xi:
\|\bm\xi'-\bm\xi\|_\infty\leq\epsilon\}$.
If
\begin{equation}
\label{eq:local_robust_approx}
\lim_{\epsilon\downarrow0}\psi_\epsilon(\bm\xi)
=\phi_\pi(\bm\xi),
\qquad \forall\bm\xi\in\Xi,
\end{equation}
then, for all sufficiently large $K$, there exists
a polyhedral policy $\pi_K$ whose
partition regions have diameter at most
$\epsilon_K:=D/\lfloor K^{1/n}\rfloor$ satisfying
\begin{equation}
\label{eq:poly_policy_bound}
0\leq\phi_{\pi_K}(\bm\xi)-\phi_\pi(\bm\xi)
\leq\psi_{\epsilon_K}(\bm\xi)-\phi_\pi(\bm\xi)+K^{-1},
\qquad \forall\bm\xi\in\Xi.
\end{equation}
Hence, $\phi_{\pi_K}$ converges pointwise
to $\phi_\pi$ as $K\to\infty$.
If $\phi_\pi$ is continuous on $\Xi$, then
the convergence is uniform, and
$\lim_{K\to\infty}\mathcal R_\Xi(\phi_{\pi_K})
=\mathcal R_\Xi(\phi_\pi)$.
In particular, \eqref{eq:local_robust_approx}
holds if $f(\bm y,\cdot)$ is upper semicontinuous
on $\Xi$ for every $\bm y\in\mathcal Y$.
\end{restatable}


\noindent\textbf{Proof.}
We first construct a family of polyhedral policies $(\pi_K)_{K\in\mathbb N}$.
By condition~\eqref{eq:local_robust_approx} and the
finiteness of $\phi_\pi$, for every $\bm\xi\in\Xi$,
there exists $r_{\bm\xi}>0$ such that
$\psi_{r_{\bm\xi}}(\bm\xi)<+\infty$.
The open balls
$\{B_{r_{\bm\xi}/2,\bm\xi}\}_{\bm\xi\in\Xi}$
cover $\Xi$, where
$B_{r,\bm\xi}:=\{\bm\xi'\in\Xi:
\|\bm\xi'-\bm\xi\|_\infty<r\}$.
By compactness of $\Xi$, there exists a finite
subcover $\{B_{r_i/2,\bm\xi_i}\}_{i\in[m]}$.
Let $\lambda>0$ be a Lebesgue number of this cover.
For sufficiently large $K$, set
$m_K:=\lfloor K^{1/n}\rfloor$ and
$\epsilon_K:=D/m_K<\lambda$.
Subdividing the enclosing hypercube into
$m_K^n\leq K$ subcubes and intersecting
them with $\Xi$ yields a polyhedral partition
whose regions have diameter at most $\epsilon_K$.
Denote the resulting partition by
$\{\Xi_k\}_{k\in\mathcal K_K}$.
By the Lebesgue number property, each nonempty
$\Xi_k$ is contained in some
$B_{r_i/2,\bm\xi_i}$, and hence in
$\bar B_{r_i,\bm\xi_i}$.
Therefore, each cell admits a finite robust
recourse decision.
For each nonempty $\Xi_k$, choose
$\bm y_k^K$$\in\mathcal Y$ satisfying
\[
\sup_{\bm\xi'\in\Xi_k}
f(\bm y_k^K,\bm\xi')
\leq
\inf_{\bm y\in\mathcal Y}
\sup_{\bm\xi'\in\Xi_k}
f(\bm y,\bm\xi')+K^{-1}.
\]
For each $\bm\xi\in\Xi$, let
$\sigma(\bm\xi):=\min\{k\in\mathcal K_K:
\bm\xi\in\Xi_k\}$ and define
$\pi_K(\bm\xi):=
\bm y_{\sigma(\bm\xi)}^K$.
Since all scenarios assigned to the same
polyhedral cell share the same response,
$\pi_K$ is a polyhedral policy.
Fix $\bm\xi\in\Xi$.
Since $\Xi_{\sigma(\bm\xi)}$ has diameter
at most $\epsilon_K$, we have
$\Xi_{\sigma(\bm\xi)}
\subseteq\bar B_{\epsilon_K,\bm\xi}$.
Consequently,
\[
\phi_\pi(\bm\xi)
\leq \phi_{\pi_K}(\bm\xi)
\leq \sup_{\bm\xi'\in\Xi_{\sigma(\bm\xi)}}
f(\bm y_{\sigma(\bm\xi)}^K,\bm\xi')
\leq \inf_{\bm y\in\mathcal Y}
\sup_{\bm\xi'\in\Xi_{\sigma(\bm\xi)}}
f(\bm y,\bm\xi')+K^{-1}
\leq \psi_{\epsilon_K}(\bm\xi)
+K^{-1},
\]
which establishes~\eqref{eq:poly_policy_bound}.
By condition~\eqref{eq:local_robust_approx},
$\phi_{\pi_K}(\bm\xi)\to\phi_\pi(\bm\xi)$
as $K\to\infty$ for every $\bm\xi\in\Xi$,
implying pointwise convergence.
Suppose $\phi_\pi$ is continuous on $\Xi$.
Fix any $\eta>0$.
By condition~\eqref{eq:local_robust_approx},
for every $\bm\xi\in\Xi$, there exist
$r_{\bm\xi}>0$ and $\bm y_{\bm\xi}\in\mathcal Y$
such that
\[
\sup_{\bm\xi'\in\bar B_{r_{\bm\xi},\bm\xi}}
f(\bm y_{\bm\xi},\bm\xi')
\leq \phi_\pi(\bm\xi)+\eta/2.
\]
By continuity of $\phi_\pi$, we may reduce
$r_{\bm\xi}$ so that
$|\phi_\pi(\bm\xi')-\phi_\pi(\bm\xi)|
<\eta/2$ for all
$\bm\xi'\in\bar B_{r_{\bm\xi},\bm\xi}$.
Thus,
\[
f(\bm y_{\bm\xi},\bm\xi')
\leq\phi_\pi(\bm\xi')+\eta,
\quad
\forall\bm\xi'\in\bar B_{r_{\bm\xi},\bm\xi}.
\]
By compactness of $\Xi$, the open cover
$\{B_{r_{\bm\xi}/2,\bm\xi}\}_{\bm\xi\in\Xi}$
admits a finite subcover
$\{B_{r_i/2,\bm\xi_i}\}_{i\in[m]}$.
Let $\delta:=\min_{i\in[m]}r_i/2>0$.
For every $\bm\xi\in\Xi$ and
$0<\epsilon<\delta$, there exists $i\in[m]$
such that
$\bar B_{\epsilon,\bm\xi}
\subseteq\bar B_{r_i,\bm\xi_i}$.
Consequently, we obtain
\[
\psi_\epsilon(\bm\xi)
\leq \sup_{\bm\xi'\in\bar B_{\epsilon,\bm\xi}}
f(\bm y_{\bm\xi_i},\bm\xi')
\leq \sup_{\bm\xi'\in\bar B_{\epsilon,\bm\xi}}
\bigl(\phi_\pi(\bm\xi')+\eta\bigr)
\leq \phi_\pi(\bm\xi)+\eta
+\sup_{\bm\xi'\in\bar B_{\epsilon,\bm\xi}}
|\phi_\pi(\bm\xi')-\phi_\pi(\bm\xi)|.
\]
Since $\phi_\pi$ is uniformly continuous on
the compact set $\Xi$, and $\eta>0$ is arbitrary,
we obtain
$\|\psi_\epsilon-\phi_\pi\|_\infty\to0$.
Together with~\eqref{eq:poly_policy_bound}
and $\epsilon_K\to0$, this implies
$\|\phi_{\pi_K}-\phi_\pi\|_\infty\to0$.
The convergence of aggregated values then
follows from Lemma~\ref{lem:riskmeasure}.
Finally, if $f(\bm y,\cdot)$ is upper semicontinuous
on $\Xi$ for every $\bm y\in\mathcal Y$,
then, for any $\bm\xi\in\Xi$,
\[
\phi_\pi(\bm\xi)\leq\psi_\epsilon(\bm\xi)
\leq\sup_{\bm\xi'\in\bar B_{\epsilon,\bm\xi}}
f(\pi(\bm\xi),\bm\xi')
\to f(\pi(\bm\xi),\bm\xi)
=\phi_\pi(\bm\xi),
\]
where the convergence follows from upper
semicontinuity at $\bm\xi$ and that
$\bm\xi\in\bar B_{\epsilon,\bm\xi}$.
Thus, condition~\eqref{eq:local_robust_approx} holds.
\pfend

\noindent\textbf{Proof of Theorem~\ref{thm:convergence}.}
We first verify condition~\eqref{eq:local_robust_approx}.
Fix $\bm\xi\in\Xi$ and let
$\bm y^\star=\pi(\bm\xi)$.
By strict feasibility, there exists
$\bar{\bm y}\in\mathcal Y_0$ such that
$h_j(\bar{\bm y},\bm\xi)<0$ for all $j\in[m]$.
For any $t\in(0,1]$, define
$\bm y_t=(1-t)\bm y^\star+t\bar{\bm y}$.
By convexity,
\[
\begin{aligned}
h_j(\bm y_t,\bm\xi)
&\leq (1-t)h_j(\bm y^\star,\bm\xi)
+t h_j(\bar{\bm y},\bm\xi)<0,
\quad j\in[m],\\
g(\bm y_t,\bm\xi)
&\leq (1-t)\phi_\pi(\bm\xi)
+t g(\bar{\bm y},\bm\xi)
\to\phi_\pi(\bm\xi)
\quad\text{as }t\downarrow0.
\end{aligned}
\]
By (i), for every sufficiently small
$t>0$, there exists $r>0$ such that
$\bm y_t$ is feasible throughout
$\bar B_{r,\bm\xi}$, and
\[
\sup_{\bm\xi'\in\bar B_{r,\bm\xi}}
g(\bm y_t,\bm\xi')
\leq\phi_\pi(\bm\xi)+\eta
\]
for any prescribed $\eta>0$.
Consequently,
$\phi_\pi(\bm\xi)\leq\psi_\epsilon(\bm\xi)
\leq\phi_\pi(\bm\xi)+\eta$
for all $0<\epsilon\leq r$.
Since $\eta>0$ is arbitrary,
condition~\eqref{eq:local_robust_approx} holds,
and pointwise convergence follows from
Lemma~\ref{lem:convergence}.

Now suppose (ii) holds. For $r\geq0$,
define the uniform variation bound
\[
\omega(r):=\sup\bigl\{
|f(\bm\xi)-f(\bm\xi')|
\mid f\in\mathcal F,\ 
\bm\xi,\bm\xi'\in\Xi,\
\|\bm\xi-\bm\xi'\|_\infty\leq r
\bigr\}.
\]
By strict feasibility, continuity in
$\bm\xi$, and compactness of $\Xi$, a finite
open-cover argument yields a uniform margin
$\delta>0$ and finitely many reference decisions
such that, for every $\bm\xi\in\Xi$, one can select
$\bar{\bm y}_{\bm\xi}\in\mathcal Y_0$ satisfying
\[
h_j(\bar{\bm y}_{\bm\xi},\bm\xi)
\leq-\delta,\qquad j\in[m].
\]
Finite-valued convexity and compactness of
$\mathcal Y_0$ give a lower bound on
$g(\cdot,\bm\xi)$ for each $\bm\xi$.
Equicontinuity and compactness of $\Xi$
make this bound uniform, so that
$\phi_\pi(\bm\xi)\geq m$ for some $m>-\infty$.
Since compactness of $\Xi$ yields finitely many
reference decisions, continuity of
$g(\bm y,\cdot)$ on $\Xi$ implies
$g(\bar{\bm y}_{\bm\xi},\bm\xi)\leq B$
for some $B<\infty$.
Thus, with $M:=B-m<\infty$,
\[
g(\bar{\bm y}_{\bm\xi},\bm\xi)
-\phi_\pi(\bm\xi)\leq M,
\qquad \forall\bm\xi\in\Xi.
\]
For sufficiently small $\epsilon>0$,
let $t=\omega(\epsilon)/\delta\leq1$ and
$\bm y_t=(1-t)\pi(\bm\xi)
+t\bar{\bm y}_{\bm\xi}$.
By convexity, for every
$\bm\xi'\in\bar B_{\epsilon,\bm\xi}$,
\[
\begin{aligned}
h_j(\bm y_t,\bm\xi')
&\leq (1-t)h_j(\pi(\bm\xi),\bm\xi)
+t h_j(\bar{\bm y}_{\bm\xi},\bm\xi)
+\omega(\epsilon)
\leq -t\delta+\omega(\epsilon)\leq0,\\
g(\bm y_t,\bm\xi')
&\leq (1-t)g(\pi(\bm\xi),\bm\xi)
+t g(\bar{\bm y}_{\bm\xi},\bm\xi)
+\omega(\epsilon)
\leq \phi_\pi(\bm\xi)+tM+\omega(\epsilon).
\end{aligned}
\]
Since $\bm y_t$ is feasible for the local robust
problem, the definition of $\psi_\epsilon$ yields
\[
0\leq\psi_\epsilon(\bm\xi)-\phi_\pi(\bm\xi)
\leq tM+\omega(\epsilon)
=\left(1+\frac{M}{\delta}\right)\omega(\epsilon).
\]
Applying Lemma~\ref{lem:convergence} with
$\epsilon_K=D/\lfloor K^{1/n}\rfloor$ yields
\[
\|\phi_{\pi_K}-\phi_\pi\|_\infty
\leq
\left(1+\frac{M}{\delta}\right)
\omega(\epsilon_K)+K^{-1}
\to0.
\]
Convergence of aggregated values follows
from the $1$-Lipschitz property of
$\mathcal R_\Xi$.

Finally, under (iii), the uniform H\"older
condition gives $\omega(r)\leq Lr^\alpha$
for some $L<\infty$. Since $\alpha/n\leq1$,
the term $K^{-1}$ is dominated by
$K^{-\alpha/n}$. Hence,
\[
\|\phi_{\pi_K}-\phi_\pi\|_\infty
=O(K^{-\alpha/n}).
\]
The aggregated values converge at the
same rate by the $1$-Lipschitz property
of $\mathcal R_\Xi$.
\pfend

\noindent\textbf{Proof of Corollary~\ref{cor:two_stage}.}
Fix an optimal first-stage solution $\bm x^\star$.
Under each regularity condition (i)--(iii),
Theorem~\ref{thm:convergence} yields a sequence
of polyhedral recourse policies $(\pi_K)_{K\in\mathbb N}$
with the corresponding convergence guarantees.
Under (ii), the convergence of aggregated values gives
\[
\begin{aligned}
0\leq z_K-z^\star
\leq c(\bm x^\star)+\mathcal R_\Xi(\phi_{\pi_K})-z^\star
=\mathcal R_\Xi(\phi_{\pi_K})
-\mathcal R_\Xi(\phi_{\bm x^\star})
\leq \|\phi_{\pi_K}-\phi_{\bm x^\star}\|_\infty
\to0.
\end{aligned}
\]
Under (iii), the same bound yields
$z_K-z^\star=O(K^{-\alpha/n})$.
\pfend

\noindent\textbf{Proof of Theorem~\ref{thm:mastersub_reform}.}
Problem~\eqref{eq:ppp_linear_kpart} can be expressed equivalently by introducing indicator variables $\lambda_k(\bm\xi)\in\{0,1\}$ that encode the assignment of each scenario $\bm\xi$ to a partition region. Specifically, $\lambda_k(\bm\xi)=1$ indicates that $\bm\xi$ belongs to region $k$, and the big-$M$ constraints enforce the corresponding objective and feasibility requirements only for the active region. This yields the following semi-infinite formulation:
$$
\begin{aligned} 
 \min_{\substack{\bm{y} \in \mathcal{Y}^{\mathcal{K}},\ z \in \mathbb{R} \\
   \bm \theta \in \Theta,\ \bm\lambda}} \ z\\
 \text{s.t.}\quad &
  \iprod{\bm{Q}(\bm{y}_k), \bm{\xi}} + \bm{q}(\bm{y}_k) \le z + M\bigl(1-\lambda_{k}(\bm\xi)\bigr), 
 \quad \forall k \in\mathcal{K},\ \forall \bm\xi \in \Xi_k^{\bm \theta},\\
 & \iprod{\bm Q_j (\bm y_{k}), \bm\xi} + \bm q_j (\bm y_{k}) \le M_j\bigl(1-\lambda_{k}(\bm\xi)\bigr), 
 \quad \forall j \in J,\ \forall k \in\mathcal{K},\ \forall \bm\xi \in \Xi_k^{\bm \theta},\\
 &\bm B_k \bigl(\bm T_{\bm \theta}\bm \xi + \bm \tau_{\bm\theta}\bigr) - \bm b_k 
 \le  M'\bigl(1-\lambda_{k}(\bm\xi)\bigr)\boldsymbol{1}, 
 \quad \forall k \in\mathcal{K},\ \forall \bm\xi \in \Xi_k^{\bm \theta},\\
 & \sum_{k \in\mathcal{K}} \lambda_{k}(\bm\xi)=1, \quad \forall \bm\xi \in \Xi,\\ 
 & {\lambda}_{k}(\bm\xi)\in \{0,1\}, \quad \forall k \in\mathcal{K},\ \forall \bm\xi \in \Xi_k^{\bm \theta}. 
\end{aligned}
$$
This is a semi-infinite mixed-integer linear program, for which a standard approach is cut generation. For a master incumbent $(\bm y,\bm\theta)$, the feasibility certificate $\max_{k\in\mathcal K, j\in J} z_{kj}^{\mathrm{fea}} \le 0$ guarantees that $(\bm y,\bm\theta)$ satisfies the original constraints in \eqref{eq:ppp_linear_kpart}. Moreover, the condition $z^{\mathrm{mst}}=\max_{k\in\mathcal K} z_k^{\mathrm{opt}}$ certifies optimality: the master relaxation yields a valid lower bound, while \eqref{eq:sub-obj} provides an upper bound once \eqref{eq:sub-feas} verifies feasibility across all regions.
\pfend

\noindent\textbf{Proof of Theorem~\ref{thm:dualization}.}
Fix $k\in\mathcal K$ and $\bm y_k\in\mathcal Y$. Constraint~\eqref{eq:ppp_linear_kpart01} is equivalent to
\[
\max_{\bm\xi\in\Xi_k^{\bm\theta}}\ \iprod{\bm Q(\bm y_k),\bm\xi}+\bm q(\bm y_k)\le z,
\]
where $\Xi_k^{\bm \theta}
=\{\bm \xi\in\mathbb R^n \mid \bm B\bm\xi\le \bm b,\ \bm B_k\bm T_{\bm\theta}\bm\xi\le \bm b_k-\bm B_k\bm\tau_{\bm\theta}\}$.
Dualizing it with dual variables $\bm\pi_k\ge \bm 0$ and $\bm\mu_k\ge \bm 0$ gives
\[
  \left\{
\begin{aligned}
  \min_{\bm \pi_k \geq 0, \bm \mu_k \geq 0} &~ \bm q(\bm y_k)+\iprod{\bm b,\bm\pi_k}+\iprod{\bm b_k-\bm B_k\bm\tau_{\bm\theta},\bm\mu_k}\\
  \text{s.t.} &~ \bm B^\top\bm\pi_k+\bm T_{\bm\theta}^\top\bm B_k^\top\bm\mu_k=\bm Q(\bm y_k).
\end{aligned}
\right\}
\le z.
\]
Setting $\bm\mu_k=\bm0$ reduces the regional dual to the dual of $\max_{\bm\xi\in\Xi}\iprod{\bm Q(\bm y_k),\bm\xi}$, which is feasible since $\Xi$ is nonempty and compact. Thus, the regional dual is always feasible, ensuring strong duality in the extended-real sense, even when $\Xi_k^{\bm\theta}$ is empty.
Similarly, applying this dualization argument to each constraint $j$ in \eqref{eq:ppp_linear_kpart02} as
\[
\max_{\bm\xi\in\Xi_k^{\bm\theta}}\ \iprod{\bm Q_j(\bm y_k),\bm\xi}+\bm q_j(\bm y_k)\le 0
\]
with dual variables $\bm\pi_{kj}\ge \bm 0$ and $\bm\mu_{kj}\ge \bm 0$ gives
\[
\bm B^\top\bm\pi_{kj}+\bm T_{\bm\theta}^\top\bm B_k^\top\bm\mu_{kj}=\bm Q_j(\bm y_k),\qquad
\bm q_j(\bm y_k)+\iprod{\bm b,\bm\pi_{kj}}+\iprod{\bm b_k-\bm B_k\bm\tau_{\bm\theta},\bm\mu_{kj}}\le 0.
\]
Combining all these constraints leads to the claimed reformulation \eqref{eq:dual-vi}. Since strong duality holds for all the dualization steps, this reformulation is exact.
\pfend

\noindent\textbf{Proof of Proposition~\ref{prop:assig_optimal}.}
Fix $\{\bm y_k\}_{k\in\mathcal K}$ and define $\mathcal K_i:=\{k\in\mathcal K: g_j(\bm y_k,\bm\xi_i)\leq 0,\ \forall j\in J\}$ and $c_{ki}:=f(\bm y_k,\bm\xi_i)$. If $\mathcal K_i=\emptyset$ for some $i$, every assignment has objective value $+\infty$, and assignment optimality is immediate. Otherwise, assigning scenario $i$ to $k\in\mathcal K_i$ contributes
\[
H_i(k):=\max\left\{
c_{ki},\max_{r\neq k}(c_{ri}-M)
\right\}
\]
to the objective of \eqref{eq:assign}. Let $k^\star\in\arg\min_{k\in\mathcal K_i}c_{ki}$. For any $k\in\mathcal K_i$, $k\neq k^\star$, setting $A:=\max_{r\notin\{k,k^\star\}}(c_{ri}-M)$ gives
\[
\begin{aligned}
H_i(k^\star)
&=\max\{c_{k^\star i},c_{ki}-M,A\}\\
&\leq\max\{c_{ki},c_{k^\star i}-M,A\}
=H_i(k),
\end{aligned}
\]
where the inequality follows from $c_{k^\star i}\leq c_{ki}$ and $M\geq 0$. Since assignments are independent across scenarios, the assignment step is optimal. For finite-step termination, let $F(\bm y,\bm\lambda)$ denote the objective of \eqref{eq:assign}, and define $v(\bm\lambda):=\min_{\bm y\in\mathcal Y^{\mathcal K}} F(\bm y,\bm\lambda)$. By exact optimization of both blocks,
\[
v(\bm\lambda^{t+1})
\leq F(\bm y^t,\bm\lambda^{t+1})
\leq F(\bm y^t,\bm\lambda^t)
=v(\bm\lambda^t).
\]
Thus, the nonincreasing objective sequence admits only finitely many strict decreases, ensuring finite termination. At termination, equality throughout the preceding inequalities establishes blockwise optimality.
\pfend

\FloatBarrier
\raggedbottom

\section{Additional Experimental Results}
\label{app:app_additional_results}

\setlength{\intextsep}{8pt}
\setlength{\textfloatsep}{8pt}
\setlength{\floatsep}{8pt}

\subsection{Parallel and Sequential Runtimes}

This appendix reports comparisons of parallel and sequential runtimes for the ALP variants on the shortest path, capital budgeting, and project management problems, as detailed in Tables~\ref{tab:results_summary_nodes_edges_time}--\ref{tab:results_summary_project_management_time}, respectively.

\begin{table}[H]
\TABLE
{\raggedright Shortest path problem: parallel and sequential runtimes.\label{tab:results_summary_nodes_edges_time}}
{\small\setlength{\tabcolsep}{4pt}\renewcommand{\arraystretch}{1.2}\centering
\begin{adjustbox}{center}
\begin{tabular}{cc*{4}{c}c*{4}{c}}
\hline
\multirow{2}{*}{Nodes (\#)} & \multirow{2}{*}{Edges (\#)}
& \multicolumn{4}{c}{Parallel Time (seconds)}
&
& \multicolumn{4}{c}{Sequential Time (seconds)} \\
\cline{3-6}\cline{8-11}
& & SVM--E & DT--E & SVM--H & DT--H & &
SVM--E & DT--E & SVM--H & DT--H \\
\hline
20  & 114  & 30.59 & 30.64 &  5.25 & \textbf{5.24} & & 30.71 & 31.10 & 9.33 & \textbf{7.33} \\
50  & 735  & 157.03 & \textbf{26.05} & 32.75 & 32.45 & & 160.75 & \textbf{27.53} & 34.92 & 34.26 \\
60  & 1062 & 152.45 & \textbf{35.75} & 97.80 & 98.15 & & 166.27 & \textbf{38.56} & 108.62 & 103.54 \\
70  & 1449 & \textbf{46.01} & 46.18 & 140.23 & 134.03 & & \textbf{46.45} & 50.17 & 144.70 & 141.10 \\
80  & 1896 & \textbf{61.28} & 61.30 & 69.69 & 74.21 & & \textbf{61.91} & 66.79 & 81.80 & 82.45 \\
90  & 2403 & 227.56 & 227.30 & \textbf{221.15} & 222.23 & & 238.59 & 234.89 & 242.92 & \textbf{232.50} \\
100 & 2970 & 253.17 & 173.17 & \textbf{108.07} & 108.81 & & 264.54 & 180.83 & 135.51 & \textbf{123.57} \\
\hline
\end{tabular}
\end{adjustbox}}
{\begingroup\fontsize{8}{11}\selectfont\emph{Notes.}
  Parallel times are measured as the maximum completion time among all parallel tasks, while the corresponding sequential times are computed as the sum of these task durations.
\par\endgroup}
\end{table}

\begin{table}[H]
\TABLE
{\raggedright Capital budgeting problem: parallel and sequential runtimes. \label{tab:results_summary_seq}}
{\small\setlength{\tabcolsep}{4pt}\renewcommand{\arraystretch}{1.2}\centering
\begin{adjustbox}{center}
\begin{tabular}{c*{4}{c}c*{4}{c}}
\hline
\multirow{2}{*}{Projects (\#)}
& \multicolumn{4}{c}{Parallel Time (seconds)}
&
& \multicolumn{4}{c}{Sequential Time (seconds)} \\
\cline{2-5}\cline{7-10}
& SVM--E & DT--E & SVM--H & DT--H & &
SVM--E & DT--E & SVM--H & DT--H \\
\hline
5    & 4.29 & 16.78 & 0.21 & \textbf{0.14} & & 4.33 & 16.78 & 1.30 & \textbf{1.22} \\
30   & 1240.56 & 1203.27 & 0.55 & \textbf{0.38} & & 1240.74 & 1203.66 & \textbf{1.14} & 1.27 \\
50   & 1220.19 & 1206.09 & 20.96 & \textbf{0.58} & & 1224.05 & 1206.70 & 22.02 & \textbf{1.86} \\
100  & 1204.71 & 1213.52 & 33.74 & \textbf{0.94} & & 1236.10 & 1214.39 & 47.14 & \textbf{2.99} \\
1000 & 1260.47 & 1253.39 & 5.97 & \textbf{5.88} & & 1262.36 & 1258.53 & \textbf{14.74} & 15.30 \\
3000 & 1304.91 & 1360.52 & 18.69 & \textbf{18.07} & & 1310.54 & 1376.43 & \textbf{35.23} & 44.05 \\
5000 & 1387.28 & 1471.68 & 30.12 & \textbf{30.07} & & 1396.08 & 1496.94 & \textbf{58.99} & 72.55 \\
\hline
\end{tabular}
\end{adjustbox}}
{\begingroup\fontsize{8}{11}\selectfont\emph{Notes.}
Parallel times are measured as the maximum completion time among all parallel tasks, while the corresponding sequential times are computed as the sum of these task durations.
\par\endgroup}
\end{table}

\begin{table}[H]
\TABLE
{\raggedright Project management problem: parallel and sequential runtimes.\label{tab:results_summary_project_management_time}}
{\small\setlength{\tabcolsep}{4pt}\renewcommand{\arraystretch}{1.2}\centering
\begin{adjustbox}{center}
\begin{tabular}{cc*{4}{c}c*{4}{c}}
\hline
\multirow{2}{*}{$m$} & \multirow{2}{*}{Projects (\#)}
& \multicolumn{4}{c}{Parallel Time (seconds)}
&
& \multicolumn{4}{c}{Sequential Time (seconds)} \\
\cline{3-6}\cline{8-11}
& & SVM--E & DT--E & SVM--H & DT--H & &
SVM--E & DT--E & SVM--H & DT--H \\
\hline
3  & 10  & 1211.21 & 1080.53 & 0.93  & \textbf{0.91} & & 1217.27 & 1082.68 & 4.87 & \textbf{3.70 }\\
8  & 25  & 1212.12 & 467.18  & 1.89  & \textbf{1.74} & & 1212.43 & 467.72 &38.83  & \textbf{ 5.89}\\
10 & 31  & 1236.96 & 1207.55 & 2.28  & \textbf{2.07} & & 1283.02 & 1208.20 & 54.49 & \textbf{6.97} \\
20 & 61  & 1273.80 & 1215.18 & 4.22  & \textbf{3.86} & & 1424.19 & 1243.74 & 183.79 & \textbf{13.03 }\\
30 & 91  & 1310.20 & 1310.09 & 6.09  & \textbf{5.71} & & 1669.82 & 1320.42 & 384.66 & \textbf{19.48 } \\
40 & 121 & 1348.77 & 1318.73 & 8.11  & \textbf{7.07} & & 1736.96 & 1333.88 & 673.56 &  \textbf{24.53}  \\
50 & 151 & 1312.53 & 1256.97 & 10.32 & \textbf{8.98} & &  1363.37 & 1279.89 & 1042.22 &  \textbf{31.77}\\
\hline
\end{tabular}
\end{adjustbox}}
{\begingroup\fontsize{8}{11}\selectfont\emph{Notes.}
  Parallel times are measured as the maximum completion time among all parallel tasks, while the corresponding sequential times are computed as the sum of these task durations.
\par\endgroup}
\end{table}

Across the three problems, sequential runtimes remain reasonably close to their parallel counterparts in many cases, indicating that the computational bottleneck lies in non-parallelized steps. For the shortest path problem, sequential runtimes are typically within $120\%$ of the parallel runtimes, with the heuristic variants (SVM--H and DT--H) remaining the fastest in most cases at under $250$ seconds. 
For capital budgeting, sequential runtimes are also relatively close to their parallel counterparts for the exact variants, while larger differences arise for some heuristic instances.
For the project management problem, the sequential runtimes of the exact variants remain relatively close to their parallel counterparts, with SVM--E and DT--E staying within approximately $130\%$ and $102\%$ of their corresponding parallel runtimes, respectively. In contrast, the heuristic variants exhibit larger sequential-to-parallel ratios, although they remain faster than the exact variants in both parallel and sequential implementations. Collectively, these results show that the robust ALP framework, particularly in its heuristic versions, scales efficiently across all three problem classes and retains its practical runtime advantage even without full parallelization.

\FloatBarrier
\subsection{Sensitivity to $K$ and Scenario Count}

\begin{figure}[H]
\FIGURE
{%
\centering
\begin{subfigure}[b]{0.48\textwidth}
  \centering
  \includegraphics[
    width=\linewidth,
    height=0.21\textheight,
    keepaspectratio
  ]{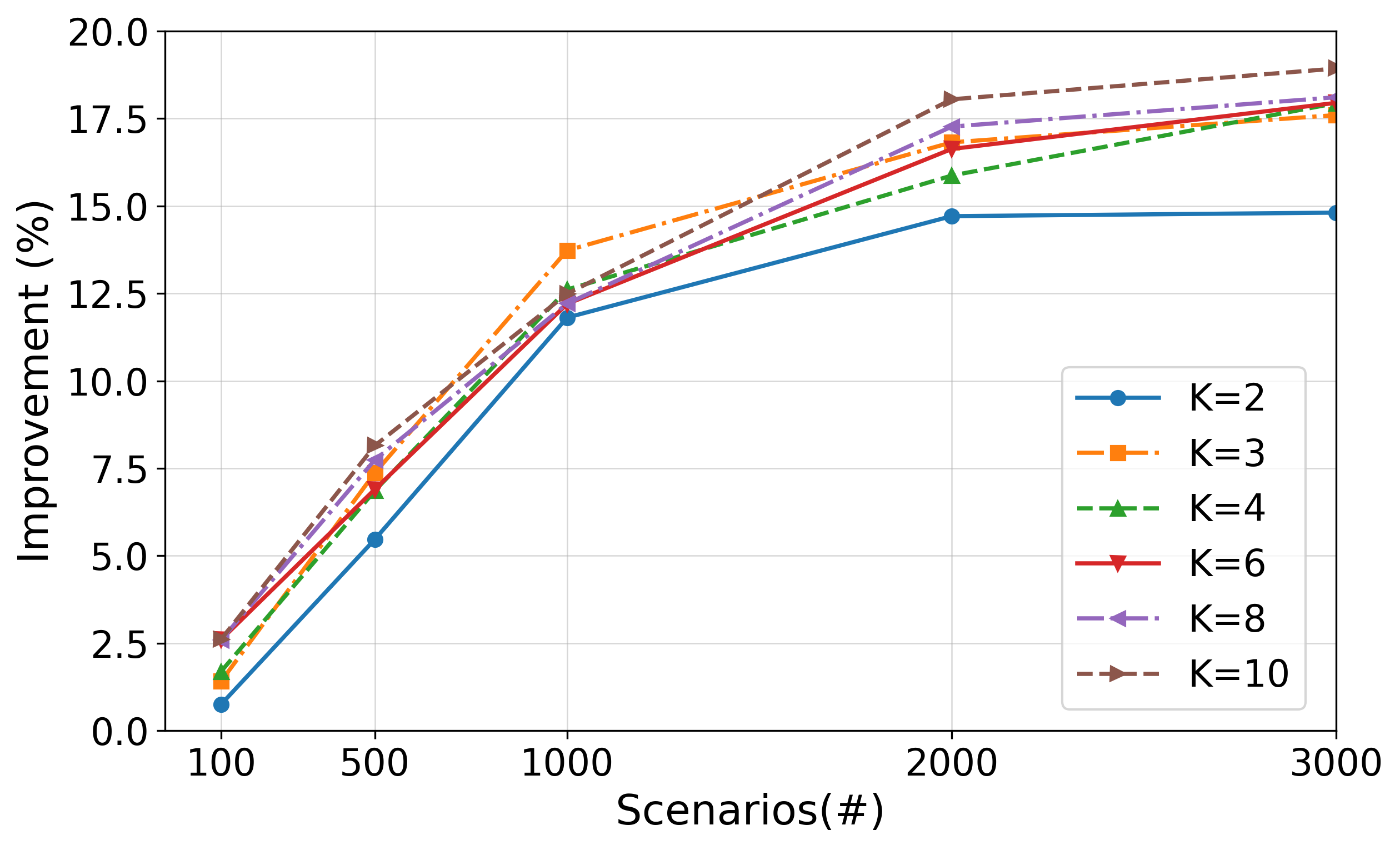}
\end{subfigure}
\hfill
\begin{subfigure}[b]{0.48\textwidth}
  \centering
  \includegraphics[
    width=\linewidth,
    height=0.21\textheight,
    keepaspectratio
  ]{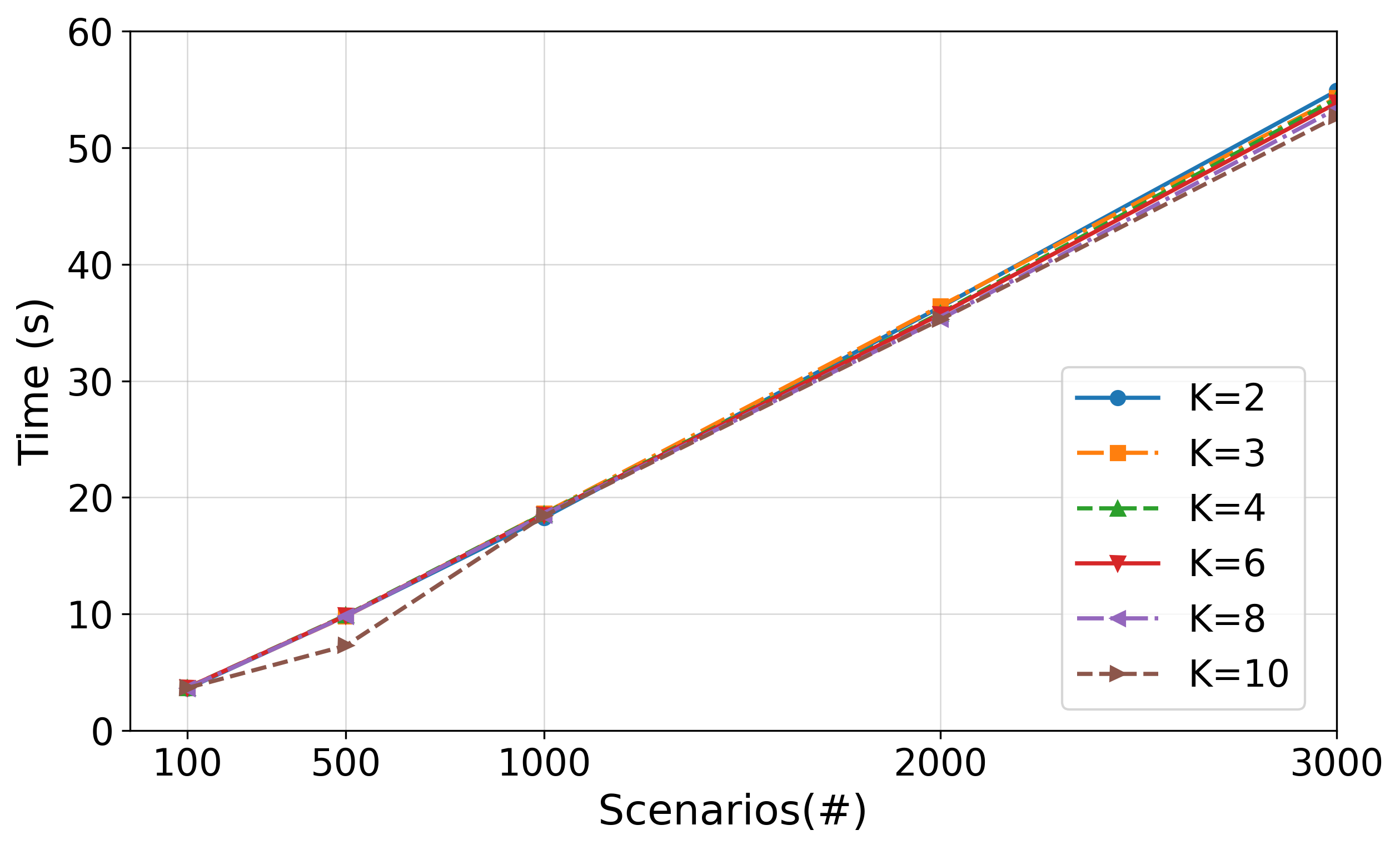}
\end{subfigure}
}
{\raggedright Sensitivity to $K$ and scenario count on the shortest path problem.\label{fig:Scenario_improvment_plot_shortest_path}}
{Experimental results of SVM--H are used for demonstration. The associated improvement percentage and runtime are averaged over all graph sizes. 
}
\end{figure}

\begin{figure}[H]
\FIGURE
{%
\centering
\begin{subfigure}[b]{0.48\textwidth}
  \centering
  \includegraphics[
    width=\linewidth,
    height=0.21\textheight,
    keepaspectratio
  ]{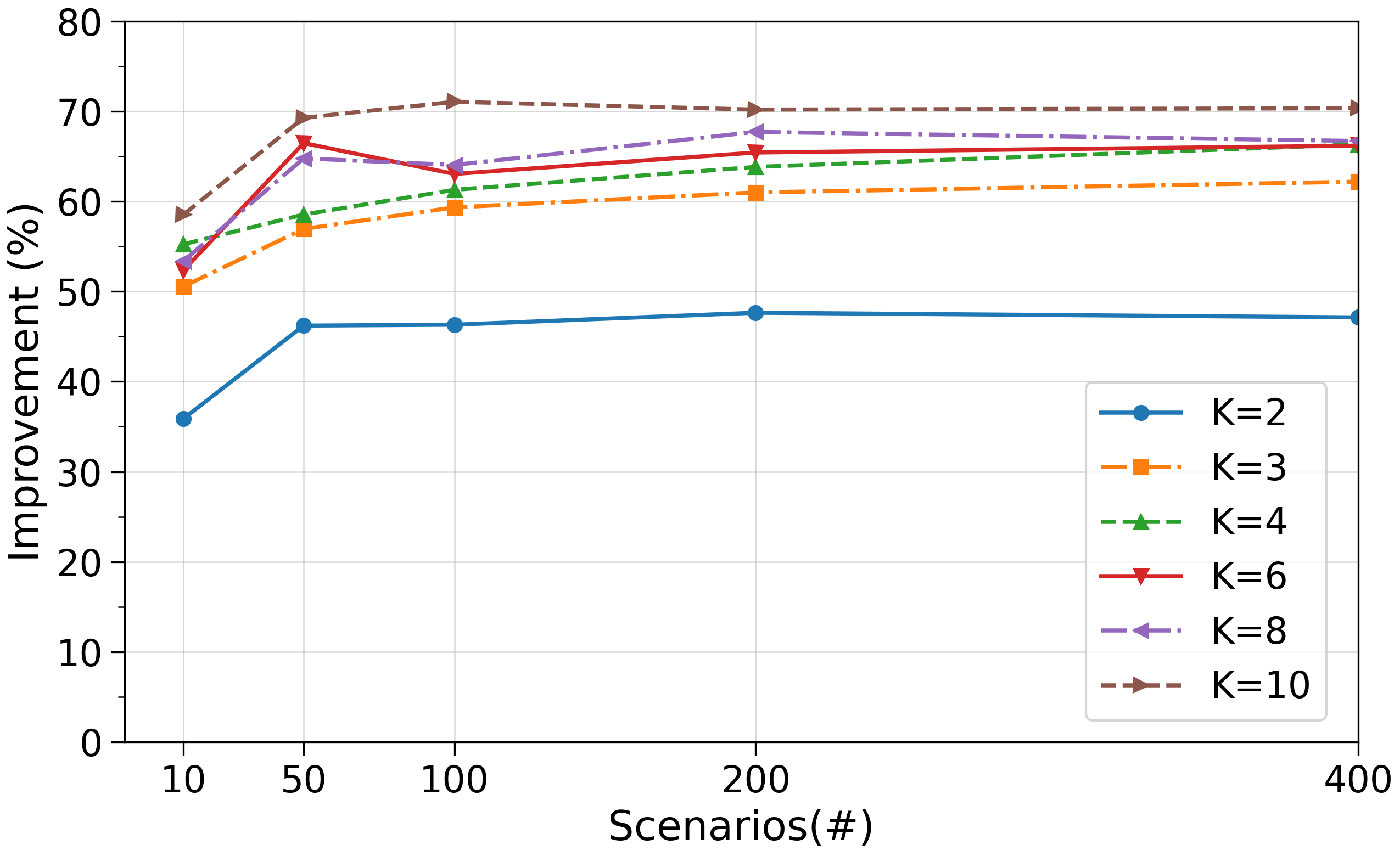}
\end{subfigure}
\hfill
\begin{subfigure}[b]{0.48\textwidth}
  \centering
  \includegraphics[
    width=\linewidth,
    height=0.21\textheight,
    keepaspectratio
  ]{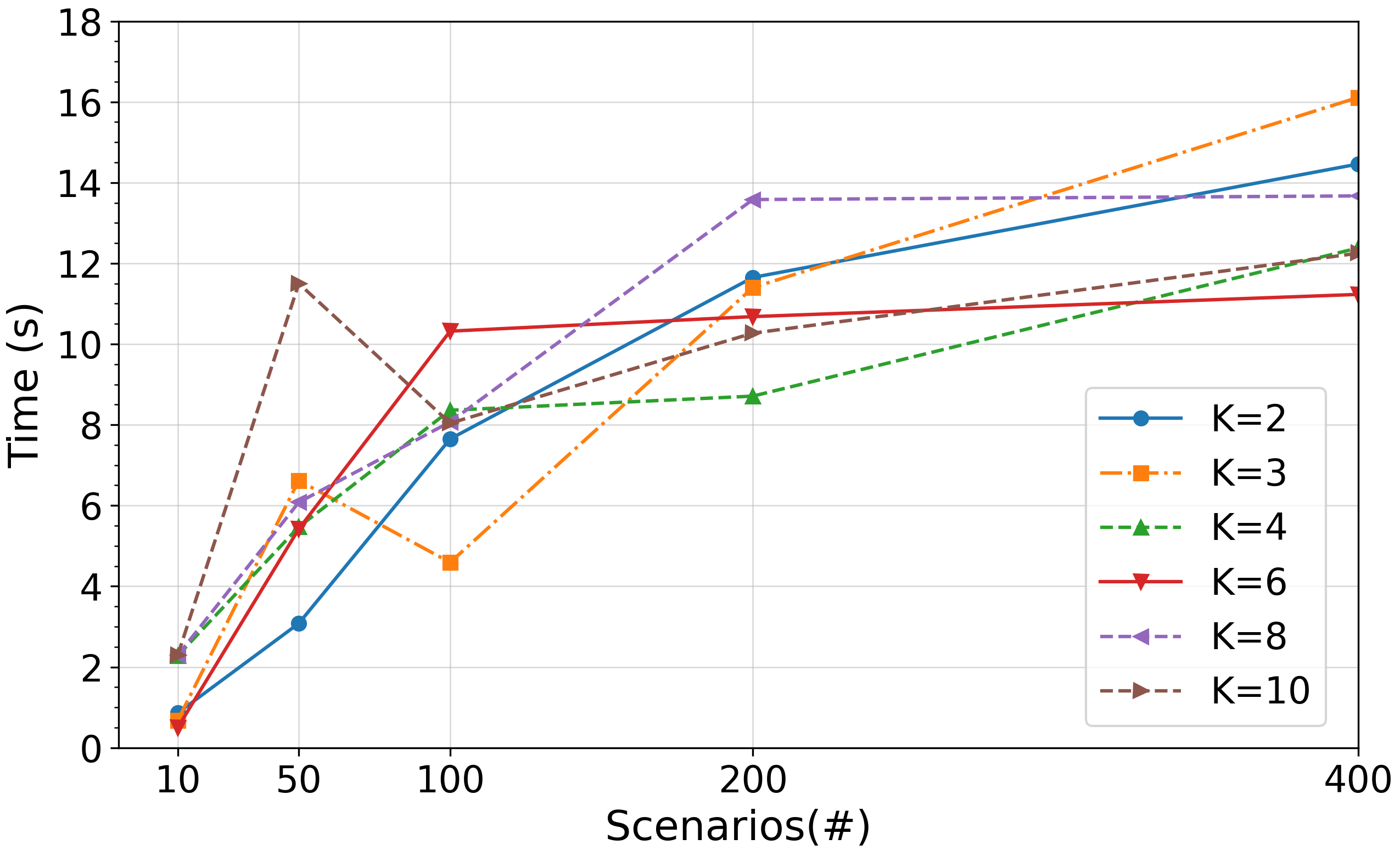}
\end{subfigure}
}
{\raggedright Sensitivity to $K$ and scenario count on the capital budgeting problem.\label{fig:Scenario_improvment_plot}}
{Experimental results of SVM--H are used for demonstration. The associated improvement percentage and runtime are averaged over all instance sizes.}
\end{figure}

\begin{figure}[H]
\FIGURE
{%
\centering
\begin{subfigure}[b]{0.48\textwidth}
  \centering
  \includegraphics[
    width=\linewidth,
    height=0.21\textheight,
    keepaspectratio
  ]{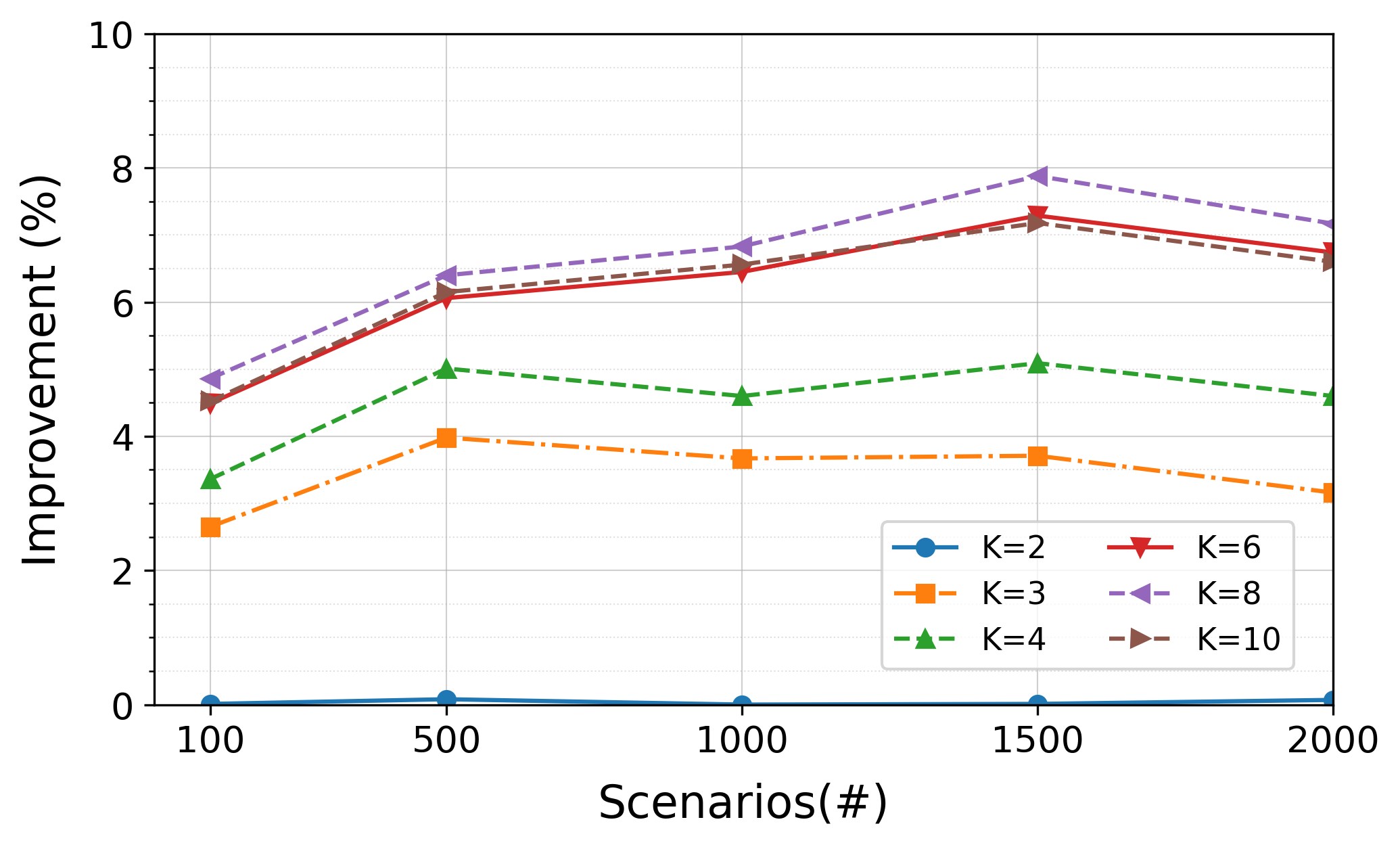}
\end{subfigure}
\hfill
\begin{subfigure}[b]{0.48\textwidth}
  \centering
  \includegraphics[
    width=\linewidth,
    height=0.21\textheight,
    keepaspectratio
  ]{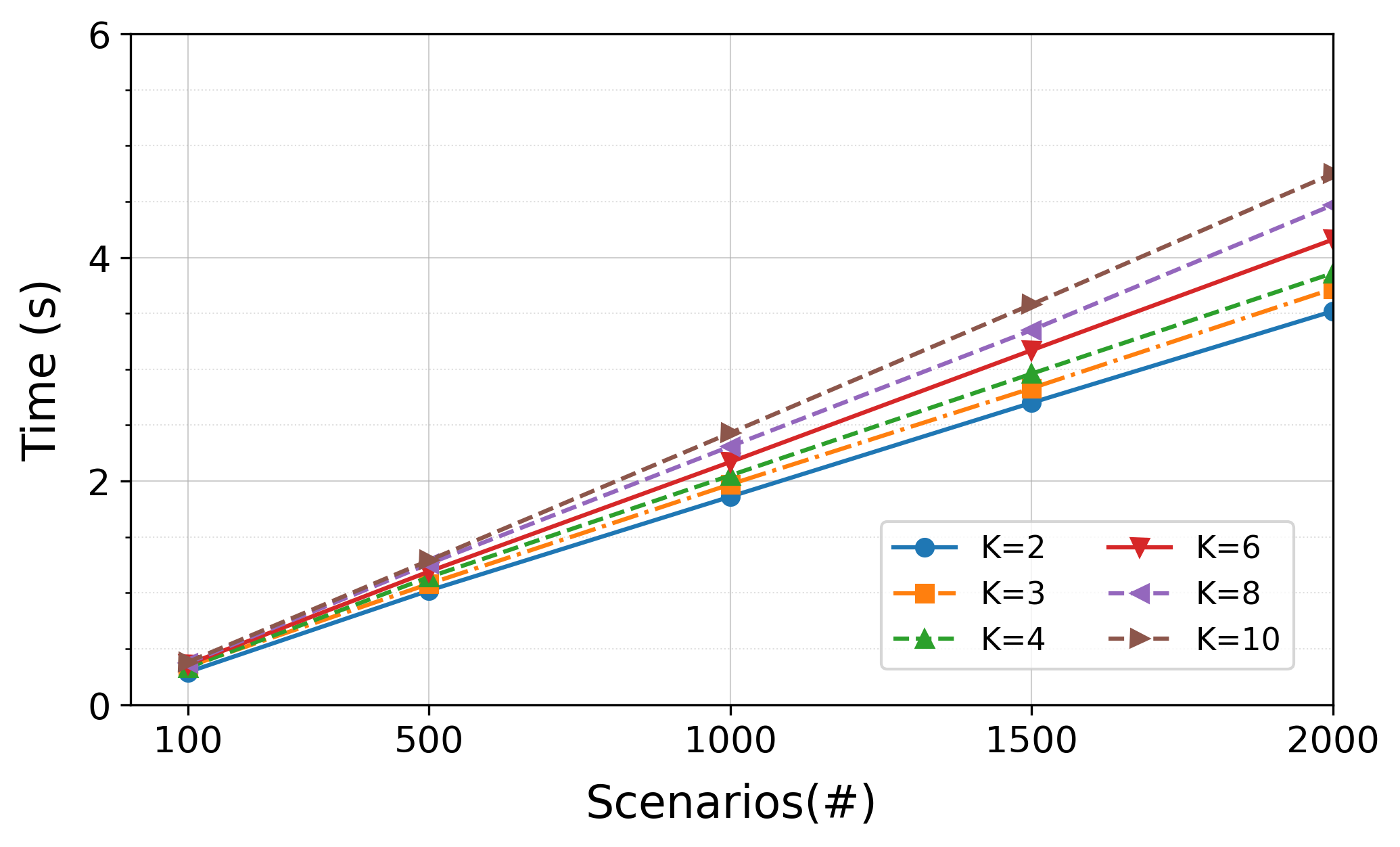}
\end{subfigure}
}
{\raggedright Sensitivity to $K$ and scenario count on the project management problem.\label{fig:Scenario_improvment_plot_project_management}}
{Experimental results of DT--H are used for demonstration. The associated improvement percentage and runtime are averaged over all instance sizes.}
\end{figure}

This subsection examines the sensitivity of ALP to the number of candidate policies $K$ and the number of scenarios used in Step~1. Results for the shortest path, capital budgeting, and project management problems are presented in Figures~\ref{fig:Scenario_improvment_plot_shortest_path}--\ref{fig:Scenario_improvment_plot_project_management}, respectively.


We use SVM--H for the first two problems and DT--H for project management, with improvement and runtime averaged over all tested instance sizes. Across all three problems, increasing $K$ and the scenario count generally improves solution quality, while runtime grows primarily with the scenario count and depends only weakly on $K$. For shortest paths, $K=10$ achieves the highest improvement in most cases. For capital budgeting, improvement increases rapidly as the scenario count grows from 10 to approximately 100, after which the gains diminish, potentially reflecting the fixed four-dimensional uncertainty set. For project management, improvements are more moderate, particularly at $K=2$, where they are nearly zero. Overall, these results suggest that larger $K$ and more scenarios can improve the learned policies, with computational cost driven primarily by the scenario count.

\end{document}